\documentclass[12pt]{article}
\usepackage{latexsym, makeidx,amssymb,amsmath}
\makeindex

\begin{document}
\baselineskip=20pt

\newcommand{\la}{\langle}
\newcommand{\ra}{\rangle}
\newcommand{\psp}{\vspace{0.4cm}}
\newcommand{\pse}{\vspace{0.2cm}}
\newcommand{\ptl}{\partial}
\newcommand{\dlt}{\delta}
\newcommand{\sgm}{\sigma}
\newcommand{\al}{\alpha}
\newcommand{\be}{\beta}
\newcommand{\G}{\Gamma}
\newcommand{\gm}{\gamma}
\newcommand{\vs}{\varsigma}
\newcommand{\Lmd}{\Lambda}
\newcommand{\lmd}{\lambda}
\newcommand{\td}{\tilde}
\newcommand{\vf}{\varphi}
\newcommand{\yt}{Y^{\nu}}
\newcommand{\wt}{\mbox{wt}\:}
\newcommand{\rd}{\mbox{Res}}
\newcommand{\ad}{\mbox{ad}}
\newcommand{\stl}{\stackrel}
\newcommand{\ol}{\overline}
\newcommand{\ul}{\underline}
\newcommand{\es}{\epsilon}
\newcommand{\dmd}{\diamond}
\newcommand{\clt}{\clubsuit}
\newcommand{\vt}{\vartheta}
\newcommand{\ves}{\varepsilon}
\newcommand{\dg}{\dagger}
\newcommand{\tr}{\mbox{Tr}}
\newcommand{\ga}{{\cal G}({\cal A})}
\newcommand{\hga}{\hat{\cal G}({\cal A})}
\newcommand{\Edo}{\mbox{End}\:}
\newcommand{\for}{\mbox{for}}
\newcommand{\kn}{\mbox{ker}}
\newcommand{\Dlt}{\Delta}
\newcommand{\rad}{\mbox{Rad}}

\begin{center}{\LARGE \bf Equivalence of Conformal Superalgebras}\end{center}
\begin{center}{\LARGE \bf to Hamiltonian Superoperators}\footnote{1991 Mathematical Subject Classification. Primary 17A 30, 17A 60; Secondary 17B 20, 81Q 60}
\end{center}
\vspace{0.2cm}

\begin{center}{\large Xiaoping Xu}\end{center}
\begin{center}{Department of Mathematics, The Hong Kong University of Science \& Technology}\end{center}
\begin{center}{Clear Water Bay, Kowloon, Hong Kong}\footnote{Research supported
 by Hong Kong RGC Competitive Earmarked Research Grant HKUST709/96P.}\end{center}

\vspace{0.3cm}

\begin{center}{\Large \bf Abstract}\end{center}
\vspace{0.2cm}

{\small In this paper, we  present  a  formal variational calculus of super functions in one real variable and find the conditions for a ``matrix differential operator'' to be a Hamiltonian superoperator. Moreover, we prove that conformal superalgebras are equivalent to certain Hamiltonian super operators. } 

\section{Introduction}

Since 1970s, Lie algebras have played more important and extensive roles in nonlinear partial differential equations and theoretical physics than they did before. One of the most interesting examples is the birth of the theory of Hamiltonian operators in middle 1970s, which was a work of Gel'fand, Dikii and Dorfman (cf. [GDi1-2], [GDo]). The existence of certain Hamiltonian operators associated with a nonlinear evolution equation implies its complete integrability. Another interesting example is the theory of vertex operator algebras introduced by Borcherds [Bo] (in initial form) and by Frenkel, Lepowsky and Meurman [FLM] (in revised form) in middle 1980s, in order to solve the problem of the moonshine representation of the Monster group. It is clear now that vertex operator algebras are the fundamental algebraic structures in conformal field theory.

  Both Hamiltonian operators and vertex operator algebras are essentially algebraic objects with one-variable structure. We observed that there should be a connection between Hamiltonian operators and vertex operator algebras many years ago. Kac [K1] introduced a concept of ``conformal superalgebra'' which is the local structure of a ``super conformal algebra'' that he and Todorov [KT] studied in middle 1980s. The theory of vertex operator superalgebras can be viewed as a restricted representation theory of the Lie algebras generated by conformal superalgebras with a Virasoro element (cf. [K1], [X3]). In this paper, we  present  a  formal variational calculus of super functions in one real variable and find the conditions for a ``matrix differential operator'' to be a Hamiltonian superoperator. Moreover, we prove that conformal superalgebras are equivalent to certain Hamiltonian super operators. 

One of the algebraic structures found in [GDo] appeared in Balinskii and Novikov's work [BN] as the local structures of certain Poisson brackets of hydrodynamic type. This is essentially a simplest example of our equivalence.

Daletsky [Da1] introduced a definition of Hamiltonian superoperator associated with an abstract complex of a Lie superalgebra. He also established in [Da1] and [Da2] a formal variational calculus over a super-commutative algebra generated by a set of so-called ``graded symbols" with coefficients in a Grassmann algebra. A deficiency of Daletsky's two works is lack of links with the other fields such as mathematical physics. In [X1], we introduced a formal variational calculus based on free fermionic fields. Furthermore, we establish in [X2] a theory of Hamiltonian superoperators of one supervariable, which is compatible with supersymmetric partial differential equations (e.g., cf. [De], [M]). Some new algebraic structures were introduced in [X2] in order to classify certain types of Hamiltonian superoperators. A connection of our Hamiltonian superoperators of one supervariable with infinite-dimensional Lie superalgebras was established. Below we shall give some technical introduction.

Throughout this paper, all the vector spaces are assumed over $\Bbb{C}$, the field of complex numbers. For two vector spaces $V_1$ and $V_2$, we denote by $LM(V_1,V_2)$ the space of linear maps from $V_1$ to $V_2$. Moreover, we denote by $\Bbb{Z}$ the ring of integers, by $\Bbb{N}$ the set of natural numbers $\{0,1,2,...\}$ and by $\Bbb{Z}_2=\Bbb{Z}/2\Bbb{Z}$ the cyclic group of order 2. When the context is clear, we use $\{0,1\}$ to denote the elements of $\Bbb{Z}_2$. We shall also use the following operator of taking residue:
$$\rd_z(z^n)=\dlt_{n,-1}\qquad\for\;\;n\in \Bbb{Z}.\eqno(1.1)$$
Furthermore, all the binomials are assumed to be expanded in the nonnegative powers of the second variable. 

 A {\it conformal superalgebra} $R=R_0\oplus R_1$ is a $\Bbb{Z}_2$-graded $\Bbb{C}[\ptl]$-module with a $\Bbb{Z}_2$-graded linear map $Y^+(\cdot,z):\;R\rightarrow LM(R,R[z^{-1}]z^{-1})$ satisfying:
$$Y^+(\ptl u,z)={dY^+(u,z)\over dz}\qquad\for\;\;u\in R;\eqno(1.2)$$
$$Y^+(u,z)v=(-1)^{ij}\rd_x{e^{x\ptl}Y^+(v,-x)u\over z-x},\eqno(1.3)$$
$$Y^+(u,z_1)Y^+(v,z_2)-(-1)^{ij}Y^+(v,z_2)Y^+(u,z_1)=\rd_x{Y^+(Y^+(u,z_1-x)v,x)\over z_2-x} \eqno(1.4)$$
for $u\in R_i;\;v\in R_j$. We denote by $(R,\ptl,Y^+(\cdot,z))$ a conformal superalgebra. When $R_1=\{0\}$, we simply call $R$ a {\it conformal algebra}. 

The above definition is the equivalent generating-function form to that given in [K1], where the author used the component formulae with $Y^+(u,z)=\sum_{n=0}^{\infty}u_{(n)}z^{-1}$.

For any two integers $m_1,m_2$,  we shall often use the following notion of index throughout this paper:
$$\ol{m_1,m_2}=\left\{\begin{array}{ll}\{m_1,m_1+1,m_1+2,...,m_2\}&\mbox{if}\;\;m_1\leq m_2,\\\emptyset&\mbox{if}\;\;m_1>m_2.\end{array}\right.\eqno(1.5)$$

 Let $({\cal G}, [\cdot,\cdot])$ be a Lie superalgebra and let $M$ be a ${\cal G}$-module. 
For a positive integer $q$,  a $q$-{\it form of} ${\cal G}$ {\it with values in} $M$ is a multi-linear map $\omega:\;{\cal G}^q={\cal G}\times \cdots \times {\cal G}\rightarrow M$ for which
$$\omega (\xi_1,\xi_2, \cdots,\xi_q)=-(-1)^{ij}\omega(\xi_1,\cdots,\xi_{\ell-1},\xi_{\ell+1},\xi_{\ell},\xi_{\ell+2},\cdots, \xi_q)\eqno(1.6)$$
for $\xi_k\in{\cal G}, \xi_{\ell}\in {\cal G}_i$ and $\xi_{\ell+1}\in {\cal G}_j.$
We denote by $c^q({\cal G},M)$ the set of $q$-forms. Moreover, we define a differential $d:\;c^q({\cal G},M)\rightarrow c^{q+1}({\cal G},M)$ by
\begin{eqnarray*} d\omega(\xi_1,\xi_2,...,\xi_{q+1})&=&\sum_{\ell=1}^{q+1}(-1)^{\ell+1+i_{\ell}(i_1+\cdots i_{\ell-1})}\xi_{\ell}\omega(\xi_1,...,\check{\xi}_{\ell},...,\xi_{q+1})\\& &+\sum_{\ell_1<\ell_2}(-1)^{\ell_1+\ell_2+(i_{\ell_1}+i_{\ell_2})(i_1+\cdots+i_{\ell_1-1})+i_{\ell_2}(i_{\ell_1+1}+\cdots+i_{\ell_2-1})}\\& &\omega([\xi_{\ell_1},\xi_{\ell_2}],\xi_1,...,\check{\xi}_{\ell_1},...,\check{\xi}_{\ell_2},....,\xi_{q+1})\hspace{4.7cm}(1.7)\end{eqnarray*}
for $\omega\in c^q({\cal G},M)$ and $\xi_k\in {\cal G}_{i_k}$ with $k\in\ol{1,q+1}$, where the above index ``check'' means deleting the term under it. A $q$-form $\omega$ is called {\it closed} if $d\omega=0$.

For any $u\in M$, we define a one-form $du$ by
$$du(\xi)=\xi(u)\qquad\for\;\;\xi\in{\cal G}.\eqno(1.8)$$
Let $\Omega$ be a subspace of $c^1({\cal G},M)$ such that $dM\subset \Omega$. Suppose that $H:\;\Omega \rightarrow {\cal G}$ is a linear map. We call $H$ $\Bbb{Z}_2$-{\it graded} if
$$H(\Omega)=H(\Omega)_0\oplus H(\Omega)_1,\qquad\mbox{where}\qquad H(\Omega)_i=H(\Omega)\bigcap{\cal G}_i.\eqno(1.9)$$
Moreover, $H$ is called {\it super skew-symmetric} if
$$\phi_1(H\phi_2)=-(-1)^{i_1i_2}\phi_2(H\phi_1)\qquad \mbox{where}\qquad H\phi_j\in H(\Omega)_{i_j}.\eqno(1.10)$$
For a $\Bbb{Z}_2$-graded super skew-symmetric linear map $H: \Omega\rightarrow {\cal G}$, we define a 2-form $\omega_H$ defined on $H(\Omega)$ by 
$$\omega_H(H\phi_1,H\phi_2)=\phi_2(H\phi_1)\qquad\for\;\;\phi_1,\phi_2\in \Omega. \eqno(1.11)$$

 We say that a super skew-symmetric $\Bbb{Z}_2$-graded linear map $H: \Omega\rightarrow {\cal G}$ is a {\it Hamiltonian superoperator}  if

(a) the subspace $H(\Omega)$ of ${\cal G}$  forms a subalgebra;

(b) the form $\omega_H$  is admissible and $d\omega_H\equiv 0$  on $H(\Omega)$.

The aim of this paper is to establish a connection between conformal superalgebras and Hamiltonian superoperators.

In Section 2, we shall establish a  formal variational calculus of super functions in one real variable and find the conditions for a ``matrix differential operator'' to be a Hamiltonian superoperator. In Section 3,  we shall  present some basic properties of conformal superalgebras. Section 4 is devoted to the proof of that  conformal superalgebras are equivalent to certain Hamiltonian super operators.

\section{Variational Calculus of Super Functions}

In this section, we shall present a  formal variational calculus of super functions in one real variable and find the conditions for a ``matrix differential operator'' to be a Hamiltonian superoperator. 

Let $\Lmd$ be a vector space that is not necessary finite-dimensional. Let $F(\Lmd)$ be the free associative algebra generated by $\Lmd$. Then the exterior algebra ${\cal E}$ generated by $\Lmd$ is isomorphic to 
$$ {\cal E}=F(\Lmd)/(\{uv+vu\mid u,v\in\Lmd\}).\eqno(2.1)$$
We can identify $\Lmd$ with its image in ${\cal E}$. Note that 
$${\cal E}={\cal E}_0\oplus  {\cal E}_1,\qquad\mbox{where}\;\;  {\cal E}_0=\sum_{n=0}^{\infty}\Lmd^{2n},\;\;{\cal E}_1=\sum_{n=0}^{\infty}\Lmd^{2n+1}.\eqno(2.2)$$
With respect to the above grading, ${\cal E}$ becomes a super-commutative associative algebra, that is,
$$uv=(-1)^{ij} vu\qquad\for\;\;u\in{\cal E}_i,\;v\in{\cal E}_j.\eqno(2.3)$$
For $i\in\Bbb{Z}_2$, let 
$$\{\psi_{i,j}\mid j\in I_i\}\eqno(2.4)$$
be a set of $C^{\infty}$-functions in a real variable $x$ with the ranges in ${\cal E}_i$, where $I_i$ is an index set. We denote 
$$\psi^{(n)}_{i,j}={d^n\psi_{i,j}\over dx^n} \qquad\for\;\;n\in\Bbb{N},\;i\in\Bbb{Z}_2,\;j\in I_i.\eqno(2.5)$$
Let ${\cal A}$ be the associative subalgebra of the algebra of functions in the real variable $x$ with the range in ${\cal E}$ generated by 
$$\{\psi^{(n)}_{i,j}\mid n\in\Bbb{N},\;i\in\Bbb{Z}_2,\;j\in I_i\}.\eqno(2.6)$$
Then ${\cal A}={\cal A}_0\oplus{\cal A}_1$ becomes a super-commutative associative algebra with
$${\cal A}_i=\mbox{span}\:\{\psi_{i_1,j_1}^{(n_1)}\psi_{i_2,j_2}^{(n_2)}\cdots \psi_{i_k,j_k}^{(n_k)}\mid k,n_{\ell}\in\Bbb{N},\;i_{\ell}\in\Bbb{Z}_2,\;j_{\ell}\in I_{i_{\ell}};\;\ell\in\ol{1,k};\;\sum_{\ell=1}^ki_{\ell}\equiv i\}.\eqno(2.7)$$

From now on, we treat $\{\psi_{i,j}^{(n)}\}$ as formal variables. Set 
$${\cal G}_i=\{\sum_{\ell\in\Bbb{Z}_2,\;j\in I_i,\;n\in\Bbb{N}}u_{\ell,j,n}\ptl_{\psi_{\ell,j}^{(n)}}\mid u_{\ell,j,n}\in {\cal A}_{i+\ell}\},\;\;{\cal G}={\cal G}_0+{\cal G}_1.\eqno(2.8)$$
Then ${\cal G}$ forms a Lie sub-superalgebra of $\mbox{Der}\:{\cal A}$. In fact, its Lie bracket is given by 
\begin{eqnarray*}& &[\ptl_1,\ptl_2]=\sum_{\ell_p\in\Bbb{Z}_2,\;j_p\in I_{\ell_p},\;n_p\in\Bbb{N};\;p=1,2}\\& &(u_{\ell_1,j_1,n_1}\ptl_{\psi_{\ell_1,j_1}^{(n_1)}}(v_{\ell_2,j_2,n_2})-(-1)^{i_1i_2}v_{\ell_1,j_1,n_1}\ptl_{\psi_{\ell_1,j_1}^{(n_1)}}(u_{\ell_2,j_2,n_2}))\ptl_{\psi_{\ell_2,j_2}^{(n_2)}}\hspace{3.3cm}(2.9)\end{eqnarray*}
for 
$$\ptl_1=\sum_{\ell\in\Bbb{Z}_2,\;j\in I_{\ell},\;n\in\Bbb{N}}u_{\ell,j,n}\ptl_{\psi_{\ell,j}^{(n)}}\in{\cal G}_{i_1},\;\;\ptl_2=\sum_{\ell\in\Bbb{Z}_2,\;j\in I_{\ell},\;n\in\Bbb{N}}v_{\ell,j,n}\ptl_{\psi_{\ell,j}^{(n)}}\in{\cal G}_{i_2}.\eqno(2.10)$$
Moreover,
$${d\over dx}=\sum_{\ell\in\Bbb{Z}_2,\;j\in I_{\ell},\;n\in\Bbb{N}}\psi_{\ell,j}^{(n+1)}\ptl_{\psi_{\ell,j}^{(n)}}.\eqno(2.11)$$
\psp

{\bf Lemma 2.1}. {\it For} 
$$\ptl=\sum_{i\in\Bbb{Z}_2,\;j\in I_i,\;n\in\Bbb{N}}u_{i,j,n}\ptl_{\psi_{i,j}^{(n)}}\in{\cal G},\eqno(2.12)$$
 $[\ptl,d/dx]=0$ {\it if and only if}
$$u_{i,j,n}=\left(d\over dx\right)^n(u_{i,j,0}),\qquad  n\in \Bbb{N},\;i\in\Bbb{Z}_2,\;j\in I_i.\eqno(2.13)$$

{\it Proof}. By (2.9) and (2.11), we have:
$$[\ptl,d/dx]=\sum_{\ell\in\Bbb{Z}_2,\;j\in I_{\ell},\;n\in\Bbb{N}}(u_{\ell,j,n+1}-du_{\ell,j,n}/dx)\ptl_{\psi_{i,j}^{(n)}}.\eqno(2.14)$$
So $[\ptl,d/dx]=0$ is equivalent to
$$u_{\ell,j,n+1}={du_{\ell,j,n}\over dx}\qquad\for\;\;\ell\in\Bbb{Z}_2,\;j\in I_{\ell},\;n\in\Bbb{N},\eqno(2.15)$$
which implies (2.13) by induction on $n.\qquad\Box$
\psp

For convenience, we shall use the notion
$$\vec{u}=\{u_{\ell,j}\mid \ell\in\Bbb{Z}_2,\; j\in I_{\ell}\}\eqno(2.16)$$
for $u_{\ell,j}\in {\cal A}$. Moreover, we define
$$\lmd \vec{u}+\mu \vec{v}=\{\lmd u_{\ell,j}+\mu v_{\ell,j}\mid \ell\in\Bbb{Z}_2,\; j\in I_{\ell}\},\eqno(2.17)$$
where $u_{\ell,j},v_{\ell,j}\in{\cal A}$ and $\lmd,\mu\in\Bbb{C}$. Set
$$\bar{\cal G}_i=\{\vec{u}=\{u_{\ell,j}\}\mid u_{\ell,j}\in {\cal A}_{\ell+i}\},\qquad\bar{\cal G}=\bar{\cal G}_0+\bar{\cal G}_1.\eqno(2.18)$$
For any $\vec{u}\in\bar{\cal G}_i$, we define
$$ \ptl_{\vec{u}}=\sum_{\ell\in\Bbb{Z}_2,\;j\in I_{\ell},\;n\in\Bbb{N}}\left({d\over dx}\right)^n(u_{\ell,j})\ptl_{\psi_{\ell,j}^{(n)}}\in {\cal G}_i.\eqno(2.19)$$
Then 
$$[\ptl_{\vec{u}},d/dx]=0\qquad\for\;\;\vec{u}\in\Bar{\cal G}\eqno(2.20)$$
by the above lemma. Moreover, 
for $\vec{u}\in\bar{\cal G}_{i_1}$ and $\vec{v}\in \bar{\cal G}_{i_2}$,
\begin{eqnarray*}& &[\ptl_{\vec{u}},\ptl_{\vec{v}}]\\&=&\sum_{\ell\in\Bbb{Z}_2,\;j\in I_{\ell},\;n\in\Bbb{N}}(\ptl_{\vec{u}}(d/dx)^n(v_{\ell,j})-(-1)^{i_1i_2}\ptl_{\vec{v}}(d/dx)^n(u_{\ell,j}))\ptl_{\psi_{\ell,j}^{(n)}}\\&=&
\sum_{\ell\in\Bbb{Z}_2,\;j\in I_{\ell},\;n\in\Bbb{N}}(d/dx)^n(\ptl_{\vec{u}}(v_{\ell,j})-(-1)^{i_1i_2}\ptl_{\vec{v}}(u_{\ell,j}))\ptl_{\psi_{\ell,j}^{(n)}}
\\&=&\ptl_{[\ptl_{\vec{u}}(\vec{v})-(-1)^{i_1i_2}\ptl_{\vec{v}}(\vec{u})]}\hspace{10.7cm}(2.21)\end{eqnarray*}
by (2.20), where
$$\ptl(\vec{w})=\{\ptl(w_{\ell,j})\mid \ell\in\Bbb{Z}_2,\; j\in I_{\ell}\}\qquad\for\;\;\ptl\in{\cal G},\;\vec{w}\in\bar{\cal G}.\eqno(2.22)$$
Thus we can define a Lie superalgebraic structure on $\bar{\cal G}$ by
$$[\vec{u},\vec{v}]=\ptl_{\vec{u}}(\vec{v})-(-1)^{i_1i_2}\ptl_{\vec{v}}(\vec{u})\qquad\for\;\;\vec{u}\in\bar{\cal G}_{i_1},\;\vec{v}\in\bar{\cal G}_{i_2}.\eqno(2.23)$$

Next we define {\it variational operators} on ${\cal A}$:
$$\dlt_{(i,j)}=\sum_{n=0}^{\infty}\left(-{d\over dx}\right)^n\circ \ptl_{\psi_{i,j}^{(n)}},\qquad \vec{\dlt}=\{\dlt_{(i,j)}\mid i\in\Bbb{Z}_2,\;j\in I_i\},\eqno(2.24)$$
where the notion $\circ$ denotes the composition of operators.
\psp

{\bf Lemma 2.2}. {\it For any} $u\in{\cal A}$,
$$\vec{\dlt}(u)=\vec{0}\Longleftrightarrow u=(d/dx)(v)+\lmd\qquad\mbox{for some}\;\;v\in{\cal A},\;\lmd\in\Bbb{C}.\eqno(2.25)$$

{\it Proof}. First we define an operator:
$$\Upsilon=\ptl_{\{\psi_{i,j}\}}=\sum_{i\in\Bbb{Z}_2,\;j\in I_i,\;n\in\Bbb{N}}\psi_{i,j}^{(n)}\ptl_{\psi_{i,j}^{(n)}}.\eqno(2.26)$$
Then $\Upsilon$ is a {\it degree operator} of polynomials in $\{\psi_{i,j}^{(n)}\}$, that is,
$$\Upsilon(\psi_{i_1,j_1}^{(n_1)}\cdots\psi_{i_k,j_k}^{(n_k)})=k\psi_{i_1,j_1}^{(n_1)}\cdots\psi_{i_k,j_k}^{(n_k)}\eqno(2.27)$$
for $k,n_{\ell}\in\Bbb{N},\;i_{\ell}\in\Bbb{Z}_2$ and $j_{\ell}\in I_{i_{\ell}}$. 

Suppose that $\vec{\dlt}(u)=\vec{0}$ for some $u\in{\cal A}$. Then
$$\dlt_{(i,j)}(u)=\sum_{n=0}^{\infty}\left(-{d\over dx}\right)^n \ptl_{\psi_{i,j}^{(n)}}(u)=0\qquad\for\;\;i\in\Bbb{Z}_2,\;j\in I_i.\eqno(2.28)$$
So
$$\ptl_{\psi_{i,j}}(u)=-\sum_{n=1}^{\infty}\left(-{d\over dx}\right)^n \ptl_{\psi_{i,j}^{(n)}}(u).\eqno(2.29)$$
Moreover, by the product rule of taking derivative, we have:
$$w_1(d/dx)(w_2)=(d/dx)(w_1w_2)-(d/dx)(w_1)w_2\qquad\;\for\;\;w_1,w_2\in{\cal A},\eqno(2.30)$$
which is equivalent to the ``integration by parts." Hence
\begin{eqnarray*}& &\psi_{i,j}\ptl_{\psi_{i,j}}(u)\\&=&-\psi_{i,j}\sum_{n=1}^{\infty}\left(-{d\over dx}\right)^n \ptl_{\psi_{i,j}^{(n)}}(u)\\&=&{d\over dx}\left[\sum_{n=1}^{\infty}\psi_{i,j}\left(-{d\over dx}\right)^{n-1} \ptl_{\psi_{i,j}^{(n)}}(u)\right]-\psi_{i,j}^{(1)}\ptl_{\psi_{i,j}^{(1)}}(u)\\& &-\sum_{n=1}^{\infty}\psi_{i,j}^{(1)}\left(-{d\over dx}\right)^n \ptl_{\psi_{i,j}^{(n+1)}}(u)\hspace{10cm}\end{eqnarray*}
\begin{eqnarray*}
&=&{d\over dx}\left[\sum_{n=1}^{\infty}\left(\psi_{i,j}\left(-{d\over dx}\right)^{n-1} \ptl_{\psi_{i,j}^{(n)}}(u)+\psi_{i,j}^{(1)}\left(-{d\over dx}\right)^{n-1} \ptl_{\psi_{i,j}^{(n+1)}}(u)\right)\right]\\& &-(\psi_{i,j}^{(1)}\ptl_{\psi_{i,j}^{(1)}}+\psi_{i,j}^{(2)}\ptl_{\psi_{i,j}^{(2)}})(u)-\sum_{n=1}^{\infty}\psi_{i,j}^{(2)}\left(-{d\over dx}\right)^n \ptl_{\psi_{i,j}^{(n+2)}}(u)\\&=&\cdots\\&=&-(\sum_{n=1}^{\infty}\psi_{i,j}^{(n)}\ptl_{\psi_{i,j}^{(n)}})(u)+(d/dx)(w_{i,j})\hspace{7.8cm}(2.31)\end{eqnarray*}
by (2.29) and (2.30), where
$$w_{i,j}=\sum_{n=1}^{\infty}\sum_{m=0}^{\infty}\psi_{i,j}^{(m)}\left(-{d\over dx}\right)^{n-1} \ptl_{\psi_{i,j}^{(n+m)}}(u).\eqno(2.32)$$
Thus
$$\Upsilon (u)=(d/dx)(w)\qquad \mbox{with}\;w=\sum_{i\in\Bbb{Z}_2,\;j\in I_i}w_{i,j}.\eqno(2.33)$$
Moreover, (2.20) and (2.26) imply
$$[\Upsilon,d/dx]=0.\eqno(2.34)$$
Set
$${\cal A}^{\dg}=\sum_{i\in\Bbb{Z}_2,\;j\in I_i,\;n\in\Bbb{N}}{\cal A}\psi_{i,j}^{(n)}.\eqno(2.35)$$
Then
$${\cal A}={\cal A}^{\dg}\oplus \Bbb{C}.\eqno(2.36)$$
We write
$$u=u'+\lmd\qquad \mbox{with}\;\;u'\in{\cal A}^{\dg},\;\lmd\in\Bbb{C}.\eqno(2.37)$$
Note that $\Upsilon$ is invertible on ${\cal A}^{\dg}$ and $w\in {\cal A}^{\dg}$ by (2.32). So
$$u'=(\Upsilon|_{{\cal A}^{\dg}})^{-1}(d/dx)(w)=(d/dx)[(\Upsilon|_{{\cal A}^{\dg}})^{-1}(w)]\eqno(2.38)$$
Therefore
$$u=u'+\lmd=(d/dx)[(\Upsilon|_{{\cal A}^{\dg}})^{-1}(w)]+\lmd,\eqno(2.39)$$
that is, the second equation in (2.25) holds.

Suppose that $u=(d/dx)(v)+\lmd$. Then
\begin{eqnarray*}& &\dlt_{(i,j)}(u)\\&=&\sum_{n=0}^{\infty}\left(-{d\over dx}\right)^n \ptl_{\psi_{i,j}^{(n)}}(u)\\&=&\sum_{n=0}^{\infty}\left(-{d\over dx}\right)^n \ptl_{\psi_{i,j}^{(n)}}(d/dx)(v)\hspace{10cm}\end{eqnarray*}
\begin{eqnarray*}
&=&-\sum_{n=0}^{\infty}\left(-{d\over dx}\right)^{n+1} \ptl_{\psi_{i,j}^{(n)}}(v)+\sum_{n=0}^{\infty}\left(-{d\over dx}\right)^n [\ptl_{\psi_{i,j}^{(n)}},d/dx](v)\\&=&-\sum_{n=0}^{\infty}\left(-{d\over dx}\right)^{n+1} \ptl_{\psi_{i,j}^{(n)}}(v)+\sum_{n=1}^{\infty}\left(-{d\over dx}\right)^n \ptl_{\psi_{i,j}^{(n-1)}}(v)\\&=&0,\hspace{13.5cm}(2.40)\end{eqnarray*}
by (2.9) and (2.11). So (2.25) holds.$\qquad\Box$
\psp

Now we let 
$$\td{\cal A}={\cal A}/(d/dx)({\cal A}).\eqno(2.41)$$
We shall use $\sim$ to denote the canonical map from ${\cal A}$ to $\td{\cal A}$. Moreover, we define an action of the Lie superalgebra $\bar{\cal G}$ on $\td{\cal A}$ by
$$\vec{u}(\td{w})=(\ptl_{\vec{u}}(w))^{\sim}\qquad\for\;\;\vec{u}\in\bar{\cal G},\;w\in{\cal A}\eqno(2.42)$$
(cf. (2.18), (2.19)). This is well defined by Lemma 2.1. Thus $\td{\cal A}$ forms a $\bar{\cal G}$-module. Note that
$$((d/dx)(w_1)w_2)^{\sim}=-(w_1(d/dx)(w_2))^{\sim}\qquad\;\for\;\;w_1,w_2\in{\cal A}\eqno(2.43)$$
by (2.30). Furthermore,
\begin{eqnarray*}\vec{u}(\td{w})&=&(\sum_{\ell\in\Bbb{Z}_2,\;j\in I_{\ell},\;n\in\Bbb{N}}\left({d\over dx}\right)^n(u_{\ell,j})\ptl_{\psi_{\ell,j}^{(n)}}(w))^{\sim}\\&=&(\sum_{\ell\in\Bbb{Z}_2,\;j\in I_{\ell},\;n\in\Bbb{N}}u_{\ell,j}\left(-{d\over dx}\right)^n\ptl_{\psi_{\ell,j}^{(n)}}(w))^{\sim}\\&=&(\sum_{\ell\in\Bbb{Z}_2,\;j\in I_{\ell}}u_{\ell,j}\dlt_{(\ell,j)}(w))^{\sim}\hspace{9cm}(2.44)\end{eqnarray*}
for $\vec{u}\in\bar{\cal G}$ and $w\in{\cal A}$ by (2.43). Set
$$\Omega=\{\vec{u}=\{u_{i,j}\}\in\bar{\cal G}\mid\mbox{only finite number of}\;u_{i,j}\;\neq 0\}.\eqno(2.45)$$ 
We identify $\Omega$ with a subspace of the space $c^1(\bar{\cal G},\td{\cal A})$ of 1-forms of $\bar{\cal G}$ with values in $\td{\cal A}$ as follows:
$$\vec{u}(\vec{v})=(\sum_{\ell\in\Bbb{Z}_2,\;j\in I_{\ell}}v_{\ell,j}u_{\ell,j})^{\sim}\qquad\for\;\;\vec{u}\in\Omega,\;\vec{v}\in \bar{\cal G}.\eqno(2.46)$$
Note that (2.44) implies
$$d(\td{w})=\vec{\dlt}(w)\in\Omega\qquad\for\;\;w\in{\cal A}\eqno(2.47)$$
(cf. (2.24)). Expression (2.25) shows that $\vec{\dlt}:\td{\cal A}\rightarrow \Omega$ is a well-defined map. In particular, $d(\td{\cal A})\subset \Omega$.

Set
$$\Omega_i=\Omega\bigcap \bar{\cal G}_i\qquad\qquad\for\;\;i\in\Bbb{Z}_2\eqno(2.48)$$
(cf. (2.18)). By (2.45),
$$\Omega=\Omega_0\oplus \Omega_1.\eqno(2.49)$$

Up to this stage, we still have a difficulty to establish the general theory of Hamiltonian superoperators for the family $(\bar{\cal G},\td{\cal A},\Omega)$, due to the difference of the parities of the components of the elements in $\Omega$ and the parity of $\bar{\cal G}$ defined in (2.18). As we indicated in the introduction, our main purpose of this paper is to give an interpretation of ``conformal superalgebras" by Hamiltonian superoperators. From this point of view, we can restrict to the smaller even family $(\bar{\cal G}_0,\td{\cal A}_0,\Omega_0)$. 

Let $H:\Omega_0\rightarrow \bar{\cal G}_0$ be a linear map as follows: for $\vec{v}\in\Omega_0$,
$$H(\vec{v})_{\ell_1,j_1}=\sum_{\ell_2\in\Bbb{Z}_2,\;j_2\in I_{\ell_2}}H^{\ell_1,j_1}_{\ell_2,j_2}v_{\ell_2,j_2},\eqno(2.50)$$
where
$$H^{\ell_1,j_1}_{\ell_2,j_2}=\sum_{n=0}^{\infty}a^{\ell_1,j_1}_{\ell_2,j_2;n}\left({d\over dx}\right)^n\eqno(2.51)$$
and
$$a^{\ell_1,j_1}_{\ell_2,j_2;n}\in{\cal A}_{\ell_1+\ell_2}\;\;\mbox{and only finite number of them are nonzero}\eqno(2.52)$$
for fixed $\ell_1,\ell_2\in\Bbb{Z}_2$ and $j_1\in I_{\ell_1},\;j_2\in I_{\ell_2}$.
Such a map $H$ is called a {\it matrix differential operator}.  For $\vec{u},\vec{v}\in\Omega_0$, we have:
\begin{eqnarray*}\vec{u}(H(\vec{v}))&=&\sum_{\ell_1\in\Bbb{Z}_2,\;j_1\in I_{\ell_1}}(H(\vec{v})_{\ell_1,j_1}u_{\ell_1,j_1})^{\sim}\\&=&\sum_{\ell_p\in\Bbb{Z}_2,\;j_p\in I_{\ell_p},\;n\in\Bbb{N};\;p=1,2}a^{\ell_1,j_1}_{\ell_2,j_2;n}((d/dx)^n(v_{\ell_2,j_2})u_{\ell_1,j_1})^{\sim}\\&=&\sum_{\ell_p\in\Bbb{Z}_2,\;j_p\in I_{\ell_p},\;n\in\Bbb{N};\;p=1,2}(-1)^{\ell_1\ell_2}a^{\ell_1,j_1}_{\ell_2,j_2;n}(u_{\ell_1,j_1}(d/dx)^n(v_{\ell_2,j_2}))^{\sim}\\&=&\sum_{\ell_p\in\Bbb{Z}_2,\;j_p\in I_{\ell_p},\;n\in\Bbb{N};\;p=1,2}(-1)^{\ell_1\ell_2}((-d/dx)^n(a^{\ell_1,j_1}_{\ell_2,j_2;n}u_{\ell_1,j_1})v_{\ell_2,j_2})^{\sim}\hspace{2cm}(2.53)\end{eqnarray*}
by (2.43). Hence the skew-symmetry of $H$:
$$\vec{u}(H(\vec{v}))=-\vec{v}(H(\vec{u}))\qquad\mbox{for any}\;\;\vec{u},\vec{v}\in\Omega_0,\eqno(2.54)$$
 is equivalent to 
$$\sum_{n=0}^{\infty}a^{\ell_1,j_1}_{\ell_2,j_2;n}\left({d\over dx}\right)^n=-(-1)^{\ell_1\ell_2}\sum_{n=0}^{\infty}\left(-{d\over dx}\right)^n
\circ a^{\ell_2,j_2}_{\ell_1,j_1;n}\eqno(2.55)$$
for $\ell_1,\ell_2\in\Bbb{Z}_2$ and $j_1\in I_{\ell_1},\;j_2\in I_{\ell_2}$, where we view $a^{\ell_1,j_1}_{\ell_2,j_2;n}$ and $a^{\ell_2,j_2}_{\ell_1,j_1;n}$ as the left multiplication operators.

Suppose that $H$ is a skew-symmetric differential operator of the form (2.50)-(2.52). We want to find the condition for $H$ to be a Hamiltonian superoperator. For $\vec{u}\in\Omega_0$, we define a linear $\ptl_{\vec{u}}(H):\bar{\cal G}_0\rightarrow \bar{\cal G}_0$ by
$$[\ptl_{\vec{u}}(H)(\vec{v})]_{\ell_1,j_1}=\sum_{\ell_2\in\Bbb{Z}_2,\;j_2\in I_{\ell_2};\;m\in\Bbb{N}}\ptl_{\vec{u}}(a^{\ell_1,j_1}_{\ell_2,j_2;m})\left({d\over dx}\right)^m(v_{\ell_2,j_2}),\eqno(2.56)$$
for $\vec{v}\in\bar{\cal G}_0,\;\ell_1\in\Bbb{Z}_2$ and $j_1\in I_{\ell_1}$.
 Now for $\vec{u},\vec{v}\in \Omega_0$, we have
\begin{eqnarray*}& &[\ptl_{H(\vec{u})}(H(\vec{v}))]_{\ell_1,j_1}\\&=&\sum_{\ell_2\in\Bbb{Z}_2,\;j_2\in I_{\ell_2},\;m\in\Bbb{N}}\ptl_{H(\vec{u})}(a^{\ell_1,j_1}_{\ell_2,j_2;m}(d/dx)^m(v_{\ell_2,j_2}))
\\&=&\sum_{\ell_2\in\Bbb{Z}_2,\;j_2\in I_{\ell_2},\;m\in\Bbb{N}}\ptl_{H(\vec{u})}(a^{\ell_1,j_1}_{\ell_2,j_2;m})(d/dx)^m(v_{\ell_2,j_2})\\& &+\sum_{\ell_2\in\Bbb{Z}_2,\;j_2\in I_{\ell_2},\;m\in\Bbb{N}}a^{\ell_1,j_1}_{\ell_2,j_2;m}\ptl_{H(\vec{u})}(d/dx)^m(v_{\ell_2,j_2})
\\&=&[\ptl_{\vec{u}}(H)(\vec{v})]_{\ell_1,j_1}+\sum_{\ell_2\in\Bbb{Z}_2,\;j_2\in I_{\ell_2},\;m\in\Bbb{N}}a^{\ell_1,j_1}_{\ell_2,j_2;m}(d/dx)^m\ptl_{H(\vec{u})}(v_{\ell_2,j_2})\\&=&[\ptl_{\vec{u}}(H)(\vec{v})]_{\ell_1,j_1}+[H\ptl_{H(\vec{u})}(\vec{v})]_{\ell_1,j_1}\hspace{8.2cm}(2.57)\end{eqnarray*}
Thus
$$\ptl_{H(\vec{u})}(H(\vec{v}))=\ptl_{\vec{u}}(H)(\vec{v})+H\ptl_{H(\vec{u})}(\vec{v})\qquad\for\;\;\vec{u},\vec{v}\in\Omega_0.\eqno(2.58)$$

Suppose that $H$ is a Hamiltonian superoperator. Let $\vec{u}_p\in\Omega_0$ with $p\in\ol{1,3}.$ Then we have
\begin{eqnarray*} H(\vec{u}_1)\omega_H(H(\vec{u})_2,H(\vec{u}_3))&=&
\ptl_{H(\vec{u}_1)}[\vec{u}_3(H(\vec{u}_2))]\\&=&\ptl_{H(\vec{u}_1)}(\vec{u}_3)(H(\vec{u}_2))+\vec{u}_3[\ptl_{H(\vec{u}_1)}(H(\vec{u}_2))]\\&=&-\vec{u}_2[H\ptl_{H(\vec{u}_1)}(\vec{u}_3)]+\vec{u}_3[\ptl_{H(\vec{u}_1)}(H(\vec{u}_2))]\hspace{2.8cm}(2.59)\end{eqnarray*}
(cf. (1.11), (2.54)) and
\begin{eqnarray*}& &\omega_H([H(\vec{u}_1),H(\vec{u}_2)],H(\vec{u}_3))\\&=&\vec{u}_3([H(\vec{u}_1),H(\vec{u}_2)])\\&=&\vec{u}_3[\ptl_{H(\vec{u}_1)}(H(\vec{u}_2))-\ptl_{H(\vec{u}_2)}(H(\vec{u}_1))]\\&=&\vec{u}_3[\ptl_{H(\vec{u}_1)}(H(\vec{u}_2))]-\vec{u}_3[\ptl_{H(\vec{u}_2)}(H)(\vec{u}_1)]-\vec{u}_3[H\ptl_{H(\vec{u}_2)}(\vec{u}_1)]\hspace{3.9cm}(2.60)\end{eqnarray*}
by Condition (a) in the abstract definition of Hamiltonian superoperators (cf. below (1.11)), (1.11), (2.23), (2.46) and (2.58). Moreover, (1.7), (1.11), (2.54), (2.59) and (2.60) imply
\begin{eqnarray*}& &d\omega_H(H(\vec{u}_1),H(\vec{u}_2),H(\vec{u}_3))\\&=&H(\vec{u}_1)\omega_H(H(\vec{u}_2),H(\vec{u}_3))+
H(\vec{u}_2)\omega_H(H(\vec{u}_3),H(\vec{u}_1))+H(\vec{u}_3)\omega_H(H(\vec{u}_1),H(\vec{u}_2))\\& &-\omega_H([H(\vec{u}_1),H(\vec{u}_2)],H(\vec{u}_3))-\omega_H([H(\vec{u}_2),H(\vec{u}_3)],H(\vec{u}_1))\\& &-\omega_H([H(\vec{u}_3),H(\vec{u}_1)],H(\vec{u}_2))\\&=&
-\vec{u}_2[H\ptl_{H(\vec{u}_1)}(\vec{u}_3)]+\vec{u}_3[\ptl_{H(\vec{u}_1)}(H(\vec{u}_2))]-\vec{u}_3[H\ptl_{H(\vec{u}_2)}(\vec{u}_1)]
\\& &+\vec{u}_1[\ptl_{H(\vec{u}_2)}(H(\vec{u}_3))]-\vec{u}_1[H\ptl_{H(\vec{u}_3)}(\vec{u}_2)]+\vec{u}_2[\ptl_{H(\vec{u}_3)}(H(\vec{u}_1))]\\& &-\vec{u}_3[\ptl_{H(\vec{u}_1)}(H(\vec{u}_2))]+\vec{u}_3[\ptl_{H(\vec{u}_2)}(H)(\vec{u}_1)]+\vec{u}_3[H\ptl_{H(\vec{u}_2)}(\vec{u}_1)] \\&&-\vec{u}_1[\ptl_{H(\vec{u}_2)}(H(\vec{u}_3))]+\vec{u}_1[\ptl_{H(\vec{u}_3)}(H)\vec{u}_2)]+\vec{u}_1[H\ptl_{H(\vec{u}_3)}(\vec{u}_2)]\\& &-\vec{u}_2[\ptl_{H(\vec{u}_3)}(H(\vec{u}_1))]+\vec{u}_2[\ptl_{H(\vec{u}_1)}(H)(\vec{u}_3)]+\vec{u}_2[H\ptl_{H(\vec{u}_1)}(\vec
{u}_3)]\\&=&\vec{u}_3[\ptl_{H(\vec{u}_2)}(H)(\vec{u}_1)]+\vec{u}_1[\ptl_{H(\vec{u}_3)}(H)(\vec{u}_2)]+\vec{u}_2[\ptl_{H(\vec{u}_1)}(H)(\vec{u}_2)]\hspace{3.5cm}(2.61)\end{eqnarray*}
for $\vec{u}_1,\vec{u}_2,\vec{u}_3\in\Omega_0$.
Hence the closedness of $\omega_H$ (cf. (1.11)) by the abstract definition of Hamiltonian superoperators below (1.11) implies
$$\vec{u}_3[\ptl_{H(\vec{u}_2)}(H)(\vec{u}_1)]+\vec{u}_1[\ptl_{H(\vec{u}_3)}(H)(\vec{u}_2)]+\vec{u}_2[\ptl_{H(\vec{u}_1)}(H)(\vec{u}_3)]=0\eqno(2.62)$$
for any $\vec{u}_1,\vec{u}_2,\vec{u}_3\in\Omega_0$.
\psp

{\bf Theorem 2.3}. {\it A matrix differential operator} $H$ {\it of the form (2.50)-(2.52) is a Hamiltonian superoperator if and only if (2.55) and (2.62) hold}.
\psp

{\it Proof.} By (2.53) and (2.61), we only need to prove that (2.55) and (2.62) imply that $H(\Omega_0)$ is a subalgebra of $\bar{\cal G}_0$. Let $\vec{u}_1,\vec{u}_2.\vec{u}_3\in \Omega_0.$ Note that
\begin{eqnarray*}& & H(\vec{u}_1)\omega_H(H(\vec{u}_2),H(\vec{u}_3))\\&=&H(\vec{u}_1)(\vec{u}_3(H(\vec{u}_2))\\&=&\ptl_{H(\vec{u}_1)}(\vec{u}_3(H(\vec{u}_2))\\&=&\ptl_{H(\vec{u}_1)}(\vec{u}_3)(H(\vec{u}_2))+\vec{u}_3(\ptl_{H(\vec{u}_1)}(H(\vec{u}_2)))\\&=&\ptl_{H(\vec{u}_1)}(\vec{u}_3)(H(\vec{u}_2))+\vec{u}_3(\ptl_{H(\vec{u}_1)}(H)(\vec{u}_2))+\vec{u}_3(H\ptl_{H(\vec{u_1})}(\vec{u_2}))\\&=&
\ptl_{H(\vec{u}_1)}(\vec{u}_3)(H(\vec{u}_2))+\vec{u}_3(\ptl_{H(\vec{u}_1)}(H)(\vec{u}_2))-\ptl_{H(\vec{u_1})}(\vec{u_2})(H(\vec{u}_3))\hspace{3.3cm}(2.63)\end{eqnarray*}
by (1.11), (2.54) and (2.58). Exchanging indices $2$ and $3$, we have
\begin{eqnarray*}& &H(\vec{u}_1)\omega_H(H(\vec{u}_3),H(\vec{u}_2))\\&=&\ptl_{H(\vec{u}_1)}(\vec{u}_2)(H(\vec{u}_3))+\vec{u}_2(\ptl_{H(\vec{u}_1)}(H)(\vec{u}_3))-\ptl_{H(\vec{u_1})}(\vec{u_3})(H(\vec{u}_2)).\hspace{3.2cm}(2.64)\end{eqnarray*}

Furthermore, 
$$H(\vec{u}_1)\omega_H(H(\vec{u}_2),H(\vec{u}_3))=-H(\vec{u}_1)\omega_H(H(\vec{u}_3),H(\vec{u}_2))\eqno(2.65)$$
by (1.11) and (2.54). Thus
\begin{eqnarray*}& &\ptl_{H(\vec{u}_1)}(\vec{u}_3)(H(\vec{u}_2))+\vec{u}_3(\ptl_{H(\vec{u}_1)}(H)(\vec{u}_2))-\ptl_{H(\vec{u_1})}(\vec{u_2})(H(\vec{u}_3))\\&=&-\ptl_{H(\vec{u}_1)}(\vec{u}_2)(H(\vec{u}_3))-\vec{u}_2(\ptl_{H(\vec{u}_1)}(H)(\vec{u}_3))+\ptl_{H(\vec{u_1})}(\vec{u_3})(H(\vec{u}_2)),\hspace{2.8cm}(2.66)\end{eqnarray*}
equivalently,
$$\vec{u}_3(\ptl_{H(\vec{u}_1)}(H)(\vec{u}_2))=-\vec{u}_2(\ptl_{H(\vec{u}_1)}(H)(\vec{u}_3)).\eqno(2.67)$$
Hence
\begin{eqnarray*}& &\vec{u}_3([H(\vec{u}_1),H(\vec{u}_2)])\\&=&\vec{u}_3(\ptl_{H(\vec{u}_1)}(H(\vec{u}_2))-\ptl_{H(\vec{u}_2)}(H(\vec{u}_1)))\\&=&\vec{u}_3[\ptl_{H(\vec{u}_1)}(H)(\vec{u}_2)+H\ptl_{H(\vec{u}_1)}(\vec{u}_2)-\ptl_{H(\vec{u}_2)}(H)(\vec{u}_1)-H\ptl_{H(\vec{u}_2)}(\vec{u}_1)]\\&=&\vec{u}_3[H(\ptl_{H(\vec{u}_1)}(\vec{u}_2)-\ptl_{H(\vec{u}_2)}(\vec{u}_1))]+\vec{u}_3[\ptl_{H(\vec{u}_1)}(H)(\vec{u}_2)]-\vec{u_3}[\ptl_{H(\vec{u}_2)}(H)(\vec{u}_1)]\\&=&\vec{u}_3[H(\ptl_{H(\vec{u}_1)}(\vec{u}_2)-\ptl_{H(\vec{u}_2)}(\vec{u}_1))]-\vec{u}_2[\ptl_{H(\vec{u}_1)}(H)(\vec{u}_3)]-\vec{u_3}[\ptl_{H(\vec{u}_2)}(H)(\vec{u}_1)]\\&=&\vec{u}_3[H(\ptl_{H(\vec{u}_1)}(\vec{u}_2)-\ptl_{H(\vec{u}_2)}(\vec{u}_1))]+\vec{u}_1[\ptl_{H(\vec{u}_3)}(H)(\vec{u}_2)]\\&=&\vec{u}_3[H(\ptl_{H(\vec{u}_1)}(\vec{u}_2)-\ptl_{H(\vec{u}_2)}(\vec{u}_1)+\vec{w})],\hspace{7.5cm}(2.68)\end{eqnarray*}
by (2.23), (2.58), (2.62) and (2.67), where we have used the fact
\begin{eqnarray*}& &\vec{u}_1[\ptl_{H(\vec{u}_3)}(H)(\vec{u}_2)]\\&=&\sum_{\ell_p\in\Bbb{Z}_2,\;j_p\in I_{\ell_p};\;p\in\ol{1,3};\;m,n\in\Bbb{N}}[(d/dx)^n(H(\vec{u}_3)_{\ell_3,j_3})\ptl_{\psi_{\ell_3,j_3}^{(n)}}(a^{\ell_1,j_1}_{\ell_2,j_2;m})(d/dx)^m(u^2_{\ell_2,j_2})u^1_{\ell_1,j_1}]^{\sim}\\&=&\sum_{\ell_p\in\Bbb{Z}_2,\;j_p\in I_{\ell_p};\;p\in\ol{1,3};\;m,n\in\Bbb{N}}\{H(\vec{u}_3)_{\ell_3,j_3}(-d/dx)^n[\ptl_{\psi_{\ell_3,j_3}^{(n)}}(a^{\ell_1,j_1}_{\ell_2,j_2;m})(d/dx)^m(u^2_{\ell_2,j_2})u^1_{\ell_1,j_1}]\}^{\sim}
\\&=&\vec{u}_3(H(\vec{w}))\hspace{12.2cm}(2.69)\end{eqnarray*}
with
$$w_{\ell_3,j_3}=-\sum_{\ell_p\in\Bbb{Z}_2,\;j_p\in I_{\ell_p};\;p=1,2;\;m,n\in\Bbb{N}}(-d/dx)^n[\ptl_{\psi_{\ell_3,j_3}^{(n)}}(a^{\ell_1,j_1}_{\ell_2,j_2;m})(d/dx)^m(u^2_{\ell_2,j_2})u^1_{\ell_1,j_1}]\eqno(2.70)$$
by (2.43), (2.54) and (2.56) with $\vec{u}_1=\{u^1_{\ell,j}\}$ and $\vec{u}_2=\{u^2_{\ell,j}\}$. Since $\vec{u}_3$ is arbitrary, we obtain
$$[H(\vec{u}_1),H(\vec{u}_2)]=H(\ptl_{H(\vec{u}_1)}(\vec{u}_2)-\ptl_{H(\vec{u}_2)}(\vec{u}_1)+\vec{w}).\eqno(2.71)$$
Therefore, $H(\Omega_0)$ forms a Lie subalgebra of $\bar{\cal G}.\qquad\Box$
\psp

We set
$${\cal L}_i=\sum_{j\in I_i}\Bbb{C}\psi_{i,j}\subset{\cal A},\qquad {\cal L}={\cal L}_0\oplus {\cal L}_1.\eqno(2.72)$$
Below we shall give two simplest examples of Hamiltonian superoperators of the form (2.50)-(2.52).
\psp

{\bf Example 2.1}. Let $k_0,k_1\in\Bbb{N}$. For $i\in\Bbb{Z}_2$, let $\{\la\cdot,\cdot\ra_{i,m}\mid m\in \ol{0,k_i}\}$ be a family of $\Bbb{C}$-bilinear forms on ${\cal L}_i$ such that 
$$\la u,v\ra_{i,m}=(-1)^{1+i+m}\la v,u\ra_{i,m}\qquad\;\;\for\;\;u,v\in{\cal L}_i.\eqno(2.73)$$
 We define the map $H:\Omega_0\rightarrow\bar{\cal G}_0$ by
$$H(\vec{u})_{i,j_1}=\sum_{m\in\ol{0,k_i},j_2\in I_i}\la\psi_{i,j_1},\psi_{i,j_2}\ra_{i,m}\left({d\over dx}\right)^m(u_{i,j_2})\eqno(2.74)$$
for $\vec{u}\in \Omega_0,\;i\in\Bbb{Z}_2,\;j_1\in I_i.$
Then $H$ is a Hamiltonian superoperator by the above theorem (cf. (2.55)).
\psp

{\bf Example 2.2}. Let $H:\Omega_0\rightarrow\bar{\cal G}_0$ be a linear map defined by (2.50) with
$$H^{\ell_1,j_1}_{\ell_2,j_2}=\sum_{j_3\in I_{\ell_1+\ell_2}}\lmd^{j_3}_{\ell_1,j_1;\ell_2,j_2}\psi_{\ell_1+\ell_2,j_3}\qquad\for\;\;\ell_p\in \Bbb{Z}_p,\;j_p\in I_{\ell_p},\eqno(2.75)$$
where $\lmd^{j_3}_{\ell_1,j_1;\ell_2,j_2}\in\Bbb{C}$. Then (2.55) is equivalent to
$$\lmd^{j_3}_{\ell_1,j_1;\ell_2,j_2}=-(-1)^{\ell_1\ell_2}\lmd^{j_3}_{\ell_2,j_2;\ell_1,j_1}\qquad\for\;\;\ell_p\in \Bbb{Z}_p,\;j_p\in I_{\ell_p},\;j_3\in I_{\ell_1+\ell_2}.\eqno(2.76)$$
Moreover,
\begin{eqnarray*}\vec{u}_3[\ptl_{H(\vec{u}_2)}(H)(\vec{u}_1)]&=&\sum_{\ell_p\in\Bbb{Z}_2,\;j_p\in I_{\ell_p};\;p\in\ol{1,3};\;j_4\in I_{\ell_1+\ell_3},\;j_5\in I_{\ell_1+\ell_2+\ell_3}}\lmd_{\ell_3,j_3;\ell_1,j_1}^{j_4}\lmd_{\ell_1+\ell_3,j_4;\ell_2,j_2}^{j_5}\\& &\psi_{\ell_1+\ell_2+\ell_3,j_5}u^2_{\ell_2,j_2}u^1_{\ell_1,j_1}u^3_{\ell_3,j_3},\hspace{6.1cm}(2.77)\end{eqnarray*}
where $\vec{u}_p=\{u^p_{i,j}\}$. Thus (2.62) is equivalent to
\begin{eqnarray*}&&0=\sum_{\ell_p\in\Bbb{Z}_2,\;j_p\in I_{\ell_p};\;p\in\ol{1,3};\;j_5\in I_{\ell_1+\ell_2+\ell_3}}\\& &[\sum_{j_4\in I_{\ell_1+\ell_3}}
\lmd_{\ell_3,j_3;\ell_1,j_1}^{j_4}\lmd_{\ell_1+\ell_3,j_4;
\ell_2,j_2}^{j_5}\psi_{\ell_1+\ell_2+\ell_3,j_5}u^2_{\ell_2,j_2}u^1_{\ell_1,j_1}u^3_{\ell_3,j_3}\hspace{4cm}\end{eqnarray*}
\begin{eqnarray*}& &+
\sum_{j_4\in I_{\ell_1+\ell_2}}\lmd_{\ell_1,j_1;\ell_2,j_2}^{j_4}\lmd_{\ell_1+\ell_2,j_4;\ell_3,j_3}^{j_5}\psi_{\ell_1+\ell_2+\ell_3,j_5}u^3_{\ell_3,j_3}u^2_{\ell_2,j_2}u^1_{\ell_1,j_1}\\& &+
\sum_{j_4\in I_{\ell_2+\ell_3}}\lmd_{\ell_2,j_2;\ell_3,j_3}^{j_4}\lmd_{\ell_2+\ell_3,j_4;\ell_1,j_1}^{j_5}\psi_{\ell_1+\ell_2+\ell_3,j_5}u^1_{\ell_1,j_1}u^3_{\ell_3,j_3}u^2_{\ell_2,j_2}],\hspace{3.9cm}(2.78)\end{eqnarray*}
equivalently
\begin{eqnarray*}&&\sum_{j_4\in I_{\ell_1+\ell_3}}\lmd_{\ell_3,j_3;\ell_1,j_1}^{j_4}\lmd_{\ell_1+\ell_3,j_4;
\ell_2,j_2}^{j_5}+(-1)^{(\ell_1+\ell_2)\ell_3}\sum_{j_4\in I_{\ell_1+\ell_2}}\lmd_{\ell_1,j_1;\ell_2,j_2}^{j_4}\lmd_{\ell_1+\ell_2,j_4;\ell_3,j_3}^{j_5}\\& &+(-1)^{(\ell_1+\ell_3)\ell_2}\sum_{j_4\in I_{\ell_2+\ell_3}}\lmd_{\ell_2,j_2;\ell_3,j_3}^{j_4}\lmd_{\ell_2+\ell_3,j_4;\ell_1,j_1}^{j_5}=0.\hspace{5.7cm}(2.79)\end{eqnarray*}

We define an algebraic operation $[\cdot,\cdot]$ on ${\cal L}$ by
$$[\psi_{\ell_1,j_1},\psi_{\ell_2,j_2}]=\sum_{j_3\in I_{\ell_1+\ell_2}}\lmd_{\ell_1,j_1;\ell_2,j_2}^{j_3}\psi_{\ell_1+\ell_2,j_3}\eqno(2.80)$$
for $\ell_p\in\Bbb{Z}_2$ and $j_p\in I_{\ell_p}$. Expression (2.76) is equivalent to the super skew-symmetry of $[\cdot,\cdot]$. Multiplying (2.79) by $\psi_{\ell_1+\ell_2+\ell_3,j_3}$ and taking summation on $j_3\in I_{ 
\ell_1+\ell_2+\ell_3}$, we obtain
\begin{eqnarray*}& &[[\psi_{\ell_3,j_3},\psi_{\ell_1,j_1}],\psi_{\ell_2,j_2}]+(-1)^{(\ell_1+\ell_2)\ell_3}[[\psi_{\ell_1,j_1},\psi_{\ell_2,j_2}],\psi_{\ell_3,j_3}]\\& &+(-1)^{(\ell_1+\ell_3)\ell_2}
[[\psi_{\ell_2,j_2},\psi_{\ell_3,j_3}],\psi_{\ell_1,j_1}]=0\hspace{7.3cm}(2.81)\end{eqnarray*}
for $\ell_p\in \Bbb{Z}_2$ and $j_p\in I_{\ell_p}$. So (2.79) implies the Jacobi identity  for $({\cal L},[\cdot,\cdot])$. Conversely, (2.81) implies (2.79) which is equivalent (2.62).

Therefore, an operator $H$ of the form (2.75) is a Hamiltonian superoperator if and only if the algebra $({\cal L},[\cdot,\cdot])$ defined by (2.72) and (2.80) forms a Lie superalgebra.
\psp

Let $H_1$ and $H_2$ be matrix differential operators of the form (2.50)-(2.52). If $\lmd H_1+\mu H_2$ is Hamiltonian for any $\lmd,\mu\in \Bbb{C}$, then we call $(H_1,H_2)$ a {\it Hamiltonian pair}. For any two matrix differential operators $H_1$ and $H_2$, we define the {\it Schouten-Nijenhuis super-bracket} $[H_1,H_2]:\: \Omega^3_0\rightarrow \td{\cal A}_0$ by
\begin{eqnarray*}& &[H_1,H_2](\vec{u}_1,\vec{u}_2,\vec{u}_3)\\&=&
\vec{u}_3[\ptl_{H_1(\vec{u}_2)}(H_2)(\vec{u}_1)]+\vec{u}_1[\ptl_{H_1(\vec{u}_3)}(H_2)(\vec{u}_2)]+\vec{u}_2[\ptl_{H_1(\vec{u}_1)}(H_2)(\vec{u}_2)]\\& &\vec{u}_3[\ptl_{H_2(\vec{u}_2)}(H_1)(\vec{u}_1)]+\vec{u}_1[\ptl_{H_2(\vec{u}_3)}(H_1)(\vec{u}_2)]+\vec{u}_2[\ptl_{H_2(\vec{u}_1)}(H_1)(\vec{u}_2)]\hspace{2.8cm}(2.82)\end{eqnarray*}
for $\vec{u}_p\in\Omega_0$. 
\psp

{\bf Corollary 2.4}. {\it Matrix differential operators} $H_1$ {\it and} $H_2$ {\it  forms a Hamiltonian pair if and only if  they satisfy (2.55) and}
$$[H_1,H_1]=0,\;\;\;[H_2,H_2]=0,\;\;\;[H_1,H_2]=0.\eqno(2.83)$$

{\bf Example 2.3}. Let $H_1$ be the Hamiltonian superoperator in Example 2.1 with $k_0=k_1=0$ and redenote
$$\la \cdot,\cdot\ra_i=\la \cdot,\cdot\ra_{i,0}\qquad\;\;i\in\Bbb{Z}_2.\eqno(2.84)$$
Let $H_2$ be the Hamiltonian superoperator in Example 2.2. So $(H_1,H_2)$ forms a Hamiltonian pair if and only if $[H_1,H_2]=0$, which is equivalent to:
\begin{eqnarray*}& &0=\sum_{\ell_p\in\Bbb{Z}_2,\;j_p\in I_{\ell_p};\;p\in\ol{1,3};\ell_1+\ell_2+\ell_3=0}\\& &[\sum_{j_4\in I_{\ell_2}}\lmd_{\ell_3,j_3;\ell_1,j_1}^{j_4}\la\psi_{\ell_2,j_4},\psi_{\ell_2,j_2}\ra_{\ell_2}u^2_{\ell_2,j_2}u^1_{\ell_1,j_1}u^3_{\ell_3,j_3}\\&&+\sum_{j_4\in I_{\ell_3}}\lmd_{\ell_1,j_1;\ell_2,j_2}^{j_4}\la\psi_{\ell_3,j_4},\psi_{\ell_3,j_3}\ra_{\ell_3}u^3_{\ell_3,j_3}
u^2_{\ell_2,j_2}u^1_{\ell_1,j_1}\\& &+\sum_{j_4\in I_{\ell_1}}\lmd_{\ell_2,j_2;\ell_3,j_3}^{j_4}\la\psi_{\ell_1,j_4},\psi_{\ell_1,j_1}\ra_{\ell_1}u^1_{\ell_1,j_1}u^3_{\ell_3,j_3}u^2_{\ell_2,j_2}\hspace{6.1cm}(2.85)\end{eqnarray*}
for $\vec{u}_p=\{u_{i,j}^p\}\in\Omega_0$. Moreover, (2.85) is equivalent to
\begin{eqnarray*}& &0=\sum_{j_4\in I_{\ell_2}}\lmd_{\ell_3,j_3;\ell_1,j_1}^{j_4}\la\psi_{\ell_2,j_4},\psi_{\ell_2,j_2}\ra_{\ell_2}+(-1)^{\ell_3}\sum_{j_4\in I_{\ell_3}}\lmd_{\ell_1,j_1;\ell_2,j_2}^{j_4}\la\psi_{\ell_3,j_4},\psi_{\ell_3,j_3}\ra_{\ell_3}\\& &+(-1)^{\ell_2}\sum_{j_4\in I_{\ell_1}}\lmd_{\ell_2,j_2;\ell_3,j_3}^{j_4}\la\psi_{\ell_1,j_4},\psi_{\ell_1,j_1}\ra_{\ell_1}
=\la [\psi_{\ell_3,j_3},\psi_{\ell_1,j_1}],\psi_{\ell_2,j_2}\ra_{\ell_2}\\& &+(-1)^{\ell_3}\la [\psi_{\ell_1,j_1},\psi_{\ell_2,j_2}],\psi_{\ell_3,j_3}\ra_{\ell_3}+(-1)^{\ell_2}\la [\psi_{\ell_2,j_2},\psi_{\ell_3,j_3}],\psi_{\ell_1,j_1}\ra_{\ell_1},\hspace{3.1cm}(2.86)\end{eqnarray*}
equivalently,
\begin{eqnarray*}& &\la [\psi_{\ell_3,j_3},\psi_{\ell_1,j_1}],\psi_{\ell_2,j_2}\ra_{\ell_2}+(-1)^{\ell_3}\la [\psi_{\ell_1,j_1},\psi_{\ell_2,j_2}],\psi_{\ell_3,j_3}\ra_{\ell_3}\\& &+(-1)^{\ell_2}\la [\psi_{\ell_2,j_2},\psi_{\ell_3,j_3}],\psi_{\ell_1,j_1}\ra_{\ell_1}=0\hspace{8.1cm}(2.87)\end{eqnarray*}
for $\ell_p\in\Bbb{Z}_2$ with $\ell_1+\ell_2+\ell_3=0$ and $j_p\in I_p$. We define $\Bbb{C}$ to be the {\it one-dimensional trivial module} of $({\cal L},[\cdot,\cdot])$, that is,
$$\xi(\mu)=0\qquad\for\;\;\xi\in{\cal L},\;\mu\in\Bbb{C}.\eqno(2.88)$$
Furthermore, we define a 2-form $\omega\in c^2({\cal L},\Bbb{C})$ by
$$\omega(u,v)=\dlt_{i_1,i_2}\la u,v\ra_{i_1}\qquad\for\;\;i_1,i_2\in\Bbb{Z}_2,\;u\in{\cal L}_{i_1},\;v\in{\cal L}_{i_2}\eqno(2.89)$$
(cf. (1.6)). The expression (2.87) is equivalent to that $\omega$ is a closed 2-form (cf. (1.7)). Thus when $(H_1,H_2)$ forms a Hamiltonian pair, we have the following one-dimensional central extension $(\bar{\cal L},[\cdot,\cdot])$ of the Lie superalgebra
$({\cal L},[\cdot,\cdot])$:
$$\bar{\cal L}_0={\cal L}_0+\Bbb{C}\subset{\cal A},\;\;\bar{\cal L}_1={\cal L}_1,\qquad \bar{\cal L}=\bar{\cal L}_0+\bar{\cal L}_1\eqno(2.90)$$
$$[u_1+\mu_1,u_2+\mu_2]=[u_1,u_2]+\omega(u_1,u_2)\qquad\for\;\;u_1,u_2\in{\cal L},\;\mu_1,\mu_2\in\Bbb{C}.\eqno(2.91)$$
\pse

{\bf Remark 2.5}. (a) In Sections 3 and 4 of [X1], we have given a theory with $I_0=\emptyset$ and define Hamiltonian superoperators for $(\bar{\cal G},\Omega,\td{\cal A})$. In particular, we are able to define another type of Hamiltonian superoperators, which we called ``odd Hamiltonian superoperators."

(b) In the theory of evolution differential equations, our functions $\{\psi_{i,j}\}$ should be in a time variable and   
a space variable. However, the structures of Hamiltonian operators associated with integrable evolution differential equations are independent of the time variable. Let $H$ be a Hamiltonian operator of the form (2.50)-(2.52). The Hamiltonian of a physical system is a special function $L$ in 
$$\{\psi_{i,j}^{(n)}(x,t)\mid i\in\Bbb{Z}_2,\;j\in I_i,\;n\in\Bbb{N}\},\eqno(2.92)$$
which is a qudratic polynomial in many interesting cases. The operator $H$ is called a {\it Hamiltonian superoperator of the system} if the system is determined by the following system of partial differential equations:
$$(\psi_{\ell_1,j_1})_t=\sum_{\ell_2\in\Bbb{Z}_2,\;j_2\in I_{\ell_2}}H^{\ell_1,j_1}_{\ell_2,j_2}\dlt_{(\ell_2,j_2)}(L)\eqno(2.93)$$
for $\ell_1\in \Bbb{Z}_2$ and $j_1\in  I_{\ell_1}$ (cf. (2.24) for $\dlt_{(\ell_2,j_2)}$). 

The well-known {\it Korteweg-de Vries (KdV) equation} is
$$ u_t=u_{xxx}+6uu_x,\eqno(2.94)$$
where $u(x,t)$ is a two-variable function. The Hamiltonian operator associated to this equation is
$$H=\ptl_x^3+4u\ptl_x+2u_x\eqno(2.95)$$
with $L=u^2/2$ (cf. [GDi1]). The operator $H$ is equivalent to 
$$H=\left({d\over dx}\right)^3+4\psi_{0,1}{d\over dx}+2\psi_{0,1}^{(1)}\eqno(2.96)$$
in terms the notions in this section when $I_1=\emptyset$ and $I_0=\{1\}$. 
 
\section{Properties of Conformal Superalgebras}
In this section, we shall  present some basic properties of conformal superalgebras.

For a vector space $V$, we set 
$$V[[z^{-1},z]]=\{\sum_{n\in\Bbb{Z}}a_nz^{-n-1}\mid a_n\in V\}.\eqno(3.1)$$
Moreover, for $g\in \Bbb{C}[[z^{-1};z]]$ and $f(z)\in V[[z^{-1};z]]$ such that the usual multiplication  $g(z)f(z)$ make sense, we have
\begin{eqnarray*}\rd_z\left(g(z){df(z)\over dz}\right)&=&\rd_z\left({d(g(z)f(z))\over dz}-{dg(z)\over dz}f(z)\right)\\&=&-\rd_z\left({dg(z)\over dz}f(z)\right).\hspace{7cm}(3.2)\end{eqnarray*}
\psp

{\bf Proposition 3.1}. {\it Let} $R=R_0\oplus R_1$ {\it be a} $\Bbb{Z}_2$-{\it graded} $\Bbb{C}[\ptl]$-{\it module and 
let} $V$ {\it be a} $\Bbb{Z}_2$-{\it graded subspace of} $R$ {\it such that}
$$R=\Bbb{F}[\ptl]V.\eqno(3.3)$$
{\it Suppose that a linear map} $Y^+(\cdot,z):\;R\rightarrow LM(R,R[z^{-1}]z^{-1})$ {\it satisfies (1.1) and (1.2) and}
$$\ptl Y^+(u,z)-Y^+(u,z)\ptl=Y^+(\ptl u,z)\qquad\mbox{\it for}\;\;u\in R.\eqno(3.4)$$
{\it Then the family} $(R,\ptl,Y^+(\cdot,z))$ {\it forms a conformal superalgebra if only if (1.3) and (1.4) acting on} $V$ {\it hold for}
$u,v\in V$.
\psp

{\it Proof}. Suppose that (1.2) and (1.3) acting on $V$ hold for $u,v\in V$. Let
$$u=f(\ptl)\bar{u},\;\;v=g(\ptl)\bar{v}\qquad\mbox{with}\;\;\bar{u}\in V_{i_1},\;\bar{v}\in V_{i_2},\;f(\ptl),g(\ptl)\in\Bbb{F}[\ptl].\eqno(3.5)$$
We shall prove (1.3) and (1.4) for such $u$ and $v$ by induction on $d(f,g)=\mbox{deg}\:f+\mbox{deg}\:g$. The statement holds when $d(f,g)=0$ by assumption. Assume that it holds for $d(f,g)\leq k$ for some $k\in\Bbb{N}$. 

Let $u, v$ be any two elements of $R$ of the form (3.5) with $d(f,g)=k$. We have:
\begin{eqnarray*}& &Y^+(\ptl u,z)v\\&=&{dY^+(\ptl u,z)v\over dz}\\&=&{d\over dz}\left((-1)^{i_1i_2}\rd_x{e^{x\ptl}Y^+(v,-x)u\over z-x}\right)\\&=&-(-1)^{i_1i_2}\rd_x{e^{x\ptl}Y^+(v,-x)u\over (z-x)^2}\\&=&
-(-1)^{i_1i_2}\rd_xe^{x\ptl}Y^+(v,-x)u\ptl_{x}\left({1\over z-x}\right)\\&=&(-1)^{i_1i_2}\rd_x{d(e^{x\ptl}Y^+(v,-x))\over dx}u\left({1\over z-x}\right)\hspace{7cm}\end{eqnarray*}
\begin{eqnarray*}
&=&(-1)^{i_1i_2}\rd_xe^{x\ptl}(\ptl Y^+(v,-x)+dY(v,-x)/dx)u\left({1\over z-x}\right)
\\&=&(-1)^{i_1i_2}\rd_xe^{x\ptl}(\ptl Y^+(v,-x)-dY(v,-x)/d(-x))u\left({1\over z-x}\right)\\&=&(-1)^{i_1i_2}\rd_xe^{x\ptl}(\ptl Y^+(v,-x)-Y(\ptl v,-x))u\left({1\over z-x}\right)\\&=&(-1)^{i_1i_2}\rd_x{e^{x\ptl}Y(\ptl v,-x)\ptl u\over z-x}\hspace{9.1cm}(3.6)\end{eqnarray*}
by (1.1), (3.2) and (3.4); 
\begin{eqnarray*}& &Y^+(u,z)\ptl v\\&=&\ptl Y^+(u,z)v-Y^+(\ptl u,z)v\\&=&(-1)^{i_1i_2}\rd_x{e^{x\ptl}\ptl Y^+(v,-x)u\over z-x}-(-1)^{i_1i_2}\rd_x{e^{x\ptl}Y^+(v,-x)\ptl u\over z-x}\\&=&(-1)^{i_1i_2}\rd_x{e^{x\ptl}(\ptl Y^+(v,-x)u-Y^+(v,-x)\ptl u)\over z-x}\\&=&(-1)^{i_1i_2}\rd_x{e^{x\ptl}Y^+(\ptl v,-x)u\over z-x}\hspace{9.1cm}(3.7)\end{eqnarray*}
by (1.1) and (3.4). Hence (1.3) holds for $d(f,g)=k+1$. So (1.3) holds for $u,v\in R$ by induction.
Moreover, acting on $V$,
\begin{eqnarray*}& &Y^+(\ptl u,z_1)Y^+(v,z_2)-(-1)^{i_1i_2}Y^+(v,z_2)Y^+(\ptl u,z_1)\\&=&\ptl_{z_1}(Y^+(u,z_1)Y^+(v,z_2)-(-1)^{i_1i_2}Y^+(v,z_2)Y^+(u,z_1))\\&=&\ptl_{z_1}\left(\rd_x{Y^+(Y^+(u,z_1-x)v,x)\over z_2-x}\right)\\&=&\rd_x{Y^+(\ptl_{z_1}(Y^+(u,z_1-x))v,x)\over z_2-x}\\&=&\rd_x{Y^+((dY^+(u,z_1-x)/d(z_1-x))v,x)\over z_2-x}
\\&=&\rd_x{Y(Y^+(\ptl u,z_1-x)v,x)\over z_2-x}\hspace{9.5cm}(3.8)\end{eqnarray*}
by (1.1) and (1.4);
\begin{eqnarray*}& &Y^+(u,z_1)Y^+(\ptl v,z_2)-(-1)^{i_1i_2}Y^+(\ptl v,z_2)Y^+(u,z_1)\\&=&\ptl_{z_2}(Y^+(u,z_1)Y^+(v,z_2)-(-1)^{i_1i_2}Y^+(v,z_2)Y^+(u,z_1))\\&=&\ptl_{z_2}\left(\rd_x{Y^+(Y^+(u,z_1-x)v,x)\over z_2-x}\right)\\&=&-\rd_x{Y^+(Y^+(u,z_1-x)v,x)\over (z_2-x)^2}\\&=&-\rd_xY^+(Y^+(u,z_1-x)v,x)\ptl_x\left({1\over z_2-x}\right)\hspace{10cm}\end{eqnarray*}
\begin{eqnarray*}&=&\rd_x\ptl_{x}(Y^+(Y^+(u,z_1-x)v,x))\left({1\over z_2-x}\right)\\&=&\rd_x{Y^+(\ptl Y^+(u,z_1-x)v,x)-Y^+( Y^+(\ptl u,z_1-x)v,x)\over z_2-x}\\&=&\rd_x{Y^+((\ptl Y^+(u,z_1-x)- Y^+(\ptl u,z_1-x))v,x)\over z_2-x}\\&=&\rd_x{Y^+(Y^+(u,z_1-x)\ptl v,x)\over z_2-x}\hspace{9cm}(3.9)\end{eqnarray*}
by (1.1), (3.2) and (3.4). Thus (1.4) acting on $V$ holds for $d(f,g)=k+1$. Therefore, (1.4) acting on $V$ holds for $u,v\in R$ by induction.

Suppose that (1.4) holds for $u\in R_i$ and $v\in R_j$ when it acts on some $w\in  R$. Then we have
\begin{eqnarray*}&&[\ptl,Y^+(u,z_1)Y^+(v,z_2)-(-1)^{ij}Y^+(v,z_2)Y^+(u,z_1)]w\\&=&
\{[\ptl, Y^+(u,z_1)]Y^+(v,z_2)-(-1)^{ij}Y^+(v,z_2)[\ptl,Y^+(u,z_1)]\\& &+Y^+(u,z_1)[\ptl,Y^+(v,z_2)]-(-1)^{ij}[\ptl,Y^+(v,z_2)]Y^+(u,z_1)\}w\\&=&(\ptl_{z_1}+\ptl_{z_2})[Y^+(u,z_1)Y^+(v,z_2)-(-1)^{ij}Y^+(v,z_2)Y^+(u,z_1)]w\\&=& \rd_x (\ptl_{z_1}+\ptl_{z_2})\left[{Y^+(Y^+(u,z_1-x)v, x)w\over z_2-x}\right]\\&=&\rd_x {\ptl_{z_1}(Y^+(Y^+(u,z_1-x)v, x))w\over z_2-x}\\& &+\rd_x Y^+(Y^+(u,z_1-x)v, x)w\ptl_{z_2}\left({1\over z_2-x}\right)\\&=&\rd_x {\ptl_{z_1}(Y^+(Y^+(u,z_1-x)v, x))w\over z_2-x}\\& &-\rd_x Y^+(Y^+(u,z_1-x)v, x)w\ptl_x\left({1\over z_2-x}\right)\\&=&\rd_x {(\ptl_{z_1}+\ptl_x)(Y^+(Y^+(u,z_1-x)v, x))w\over z_2-x}\\&=&\rd_x{Y^+(\ptl Y^+(u,z_1-x)v, x)w\over z_2-x}\\&=&\rd_x{[\ptl,Y^+(Y^+(u,z_1-x)v, x)]w\over z_2-x}\\&=&\left[\ptl,\rd_x{Y^+(Y^+(u,z_1-x)v, x)\over z_2-x}\right]w\hspace{7.8cm}(3.10)\end{eqnarray*}
by (1.1), (3.2), and (3.4). Thus (1.4) holds when it acts on $\ptl w$. So it holds when it acts on the $\Bbb{F}[\ptl]$-submodule generated by $V$, which is equal to $R$ by (3.3). Therefore, (1.4) holds and $(R,\ptl,Y^+(\cdot,z))$ forms a conformal superalgebra. 

Conversely if $(R,\ptl,Y^+(\cdot,z))$ forms a conformal superalgebra, obviously (1.3) and (1.4) acting on $V$ hold for
$u,v\in V.\qquad\Box$
\psp

{\bf Remark 3.2}. Kac [K1] proved the above proposition by an indirect and relatively complicated method.
\psp

{\bf Proposition 3.3}. {\it Let} $R$ {\it be a} $\Bbb{Z}_2${\it -graded space and let} $Y^+(\cdot,z):R\rightarrow LM(R,R[z^{-1}]z^{-1})$ {\it be a linear map satisfying  (1.1). Then Expressions (1.2) and (1.3) imply (3.4)}.
\psp

{\it Proof}. For $u\in R_{i_1}$ and $v\in R_{i_2}$, we have
\begin{eqnarray*}& &(-1)^{i_1i_2}\rd_x{e^{x\ptl}Y^+(v,-x)\ptl u\over z-x}\\&=&Y^+(\ptl u,z)v\\&=&{d\over dz}Y^+(u,z)v\\&=&(-1)^{i_1i_2}{d\over dz}\rd_x{e^{x\ptl}Y^+(v,-x)u\over z-x}\\&=&(-1)^{i_1i_2}\rd_x
\ptl_z{e^{x\ptl}Y^+(v,-x)u\over z-x}\\&=&-(-1)^{i_1i_2}\rd_x{e^{x\ptl}Y^+(v,-x)u\over (z-x)^2}
\\&=&-(-1)^{i_1i_2}\rd_x\ptl_x\left({1\over z-x}\right)e^{x\ptl}Y^+(v,-x)u\\&=&(-1)^{i_1i_2}\rd_x{1\over z-x}\ptl_x(e^{x\ptl}Y^+(v,-x)u)\\&=&(-1)^{i_1i_2}\rd_x{1\over z-x}[e^{x\ptl}(\ptl Y^+(v,-x)-dY^+(v,-x)/d(-x))u\\&=&
(-1)^{i_1i_2}\rd_x{e^{x\ptl}(\ptl Y^+(v,-x)-Y^+(\ptl v,-x))u\over z-x}\hspace{5.8cm}(3.11)\end{eqnarray*}
by (3.2), (1.1) and (1.3). Hence
$$\rd_x{e^{x\ptl}[\ptl Y^+(v,-x)-Y^+(v,-x)\ptl-Y^+(\ptl v,-x) ]u\over z-x}=0.\eqno(3.12)$$
Since
$$[\ptl Y^+(v,-x)-Y^+(v,-x)\ptl-Y^+(\ptl v,-x) ]u\in R[x^{-1}]x^{-1},\eqno(3.13)$$
we can prove
$$[\ptl Y^+(v,-x)-Y^+(v,-x)\ptl -Y^+(\ptl v,-x)]u=0\eqno(3.14)$$
by induction on the degree of $z$ starting with the minimal degree. Thus we have 
$$\ptl Y^+(v,-x)-Y^+(v,-x)\ptl =Y^+(\ptl v,-x)\eqno(3.15)$$
because $u$ is arbitrary. So (3.4) holds.$\qquad\Box$
\psp

{\bf Proposition 3.4}. {\it Let} $R$ {\it be a} $\Bbb{Z}_2${\it -graded space and let} $Y^+(\cdot,z):R\rightarrow LM(R,R[z^{-1}]z^{-1})$ {\it be a linear map satisfying  (1.1)-(1.3). Then Expressions (1.4) is symmetric with respect to} $(u,z_1)$ {\it and} $(v,z_2)$.
\psp

{\it Proof}. It is enough to prove that the right-hand side of (1.4) is supersymmetric with respect to $(u,z_1)$ and $(v,z_2)$. For $u\in R_{i_1}$ and $v\in R_{i_2}$, we note that
\begin{eqnarray*}& &\rd_x{Y^+(Y^+(u,z_1-x)v,x)\over z_2-x}\\&=&(-1)^{i_1i_2}\rd_{x,y}{Y^+(e^{y\ptl}Y^+(v,-y)u,x)\over (z_2-x)(z_1-x-y)}\\&=&(-1)^{i_1i_2}\rd_{x,y}{e^{y\ptl_x}(Y^+(Y^+(v,-y)u,x))\over (z_2-x)(z_1-x-y)}\\&=&(-1)^{i_1i_2}\rd_{x,y}Y^+(Y^+(v,-y)u,x)e^{-y\ptl_x}\left({1\over (z_2-x)(z_1-x-y)}\right)\\&=&(-1)^{i_1i_2}\rd_{x,y}Y^+(Y^+(v,-y)u,x){1\over (z_2-x+y)(z_1-x)}\\&=&-(-1)^{i_1i_2}\rd_{x,y}Y^+(Y^+(v,-y)u,x){1\over (-z_2+x-y)(z_1-x)}\\&=&-(-1)^{i_1i_2}\rd_x{Y^+(Y^+(v,z_2-x)u,x)\over(z_1-x)}\hspace{7.5cm}(3.16)\end{eqnarray*}
by (1.1)-(1.3), (3.2) and (3.4).$\qquad\Box$ 

\section{Equivalence Theorem}

In this section, we shall prove that conformal superalgebras are equivalent to certain linear Hamiltonian superoperators.

 Let $R$ be a $\Bbb{Z}_2$-graded free $\Bbb{C}[\ptl]$-module over its $\Bbb{Z}_2$-graded subspace $V$ and let $Y^+(\cdot,z):V\rightarrow LM(V,R[z^{-1}]z^{-1})$ be a linear map such that  
$$Y^+(V_{i_1},z)V_{i_2}\subset R_{i_1+i_2}[z^{-1}]z^{-1}\qquad\for\;\;i_1,i_2\in\Bbb{Z}_2.\eqno(4.1)$$
We can first extend $Y^+(\cdot,z)$ to a linear map $Y^+(\cdot,z):R\rightarrow LM(V,R[z^{-1}]z^{-1})$ by
$$Y^+(f(\ptl)u,z)v=f(d/dz)Y^+(u,z)v\qquad\;\;\for\;\;u,v\in V.\eqno(4.2)$$
Then we extend $Y^+(\cdot,z)$ to a linear map $Y^+(\cdot,z):R\rightarrow LM(R,R[z^{-1}]z^{-1})$ by
$$Y^+(u,z)\ptl^mv=\sum_{p=0}^m(^m_{\;p})\ptl^{m-p}Y^+((-\ptl)^pu,z)v\qquad\for\;\;u\in R,\;v\in V,\;m\in\Bbb{N}.\eqno(4.3)$$
Note the map $Y^+(\cdot,z)$ naturally satisfies (1.1), (1.2) and (3.4). Thus if $Y^+(\cdot,z)$ satisfies (1.3) and (1.4) acting on $V$ for $u,v\in V$, then $(R,Y^+(\cdot,z))$ forms a conformal superalgebra with the extended map $Y^+(\cdot,z)$ by Proposition 3.1

Let $\{\vs_{i,j}\in I_i\}$ be a fixed basis of $V_i$ for $i\in \Bbb{Z}_2$, where $I_1$ and $I_2$ are index sets. We write
$$Y^+(\vs_{\ell_1,j_1},z)\vs_{\ell_2,j_2}=\sum_{j_3\in I_{\ell_1+\ell_2},\;m,n\in\Bbb{N}}\lmd_{\ell_1,j_1;\ell_2,j_2}^{j_3,n,m}\ptl^n\vs_{\ell_1+\ell_2,j_3}z^{-m-1}\eqno(4.4)$$
for $\ell_1,\ell_2\in \Bbb{Z}_2$ and $j_p\in I_{\ell_p}$, where $\lmd_{\ell_1,j_1;\ell_2,j_2}^{j_3,n,m}\in\Bbb{C}$.
Now (1.3) is equivalent to
\begin{eqnarray*}&&\sum_{j_3\in I_{\ell_1+\ell_2},\;m,n\in\Bbb{N}}\lmd_{\ell_1,j_1;\ell_2,j_2}^{j_3,n,m}\ptl^n\vs_{\ell_1+\ell_2,j_3}z^{-m-1}\\&=&Y^+(\vs_{\ell_1,j_1},z)\vs_{\ell_2,j_2}\\&=&(-1)^{\ell_1\ell_2}\rd_x{e^{x\ptl}Y^+(\vs_{\ell_2,j_2},-x)\vs_{\ell_1,j_1}\over z-x}\\&=&
(-1)^{\ell_1\ell_2}\rd_x{1\over z-x}(\sum_{j_3\in I_{\ell_1+\ell_2},\;m,n\in\Bbb{N}}\lmd_{\ell_2,j_2;\ell_1,j_1}^{j_3,n,m}e^{x\ptl}\ptl^n\vs_{\ell_1+\ell_2,j_3}(-x)^{-m-1})\hspace{4cm}\end{eqnarray*}
\begin{eqnarray*}&=&
(-1)^{\ell_1\ell_2}\sum_{j_3\in I_{\ell_1+\ell_2},\;m,n,p\in\Bbb{N}}{(-1)^{m+p}\over p!}\lmd_{\ell_2,j_2;\ell_1,j_1}^{j_3,n,m+p}\ptl^{n+p}\vs_{\ell_1+\ell_2,j_3}z^{-m-1}\\&=&-(-1)^{\ell_1\ell_2}\sum_{j_3\in I_{\ell_1+\ell_2},\;m,n,p\in\Bbb{N}}{(-1)^{m+p}\over p!}\lmd_{\ell_2,j_2;\ell_1,j_1}^{j_3,n-p,m+p}\ptl^n\vs_{\ell_1+\ell_2,j_3}z^{-m-1}\\&=&-(-1)^{\ell_1\ell_2}\sum_{j_3\in I_{\ell_1+\ell_2},\;m,n,p\in\Bbb{N},\;p\geq m}{(-1)^p\over (p-m)!}\lmd_{\ell_2,j_2;\ell_1,j_1}^{j_3,m+n-p,p}\ptl^n\vs_{\ell_1+\ell_2,j_3}z^{-m-1}
\hspace{2.4cm}(4.5)\end{eqnarray*}
for $\ell_p\in \Bbb{Z}_2$ and $j_p\in I_p$, where we treat
$$\lmd_{\ell_1,j_1;\ell_2,j_2}^{j_3,n,m}=0\;\;\mbox{if}\;\;m<0\;\;\mbox{or}\;\;n<0.\eqno(4.6)$$
Expression (4.5) shows that (1.3) is equivalent to
$$\lmd_{\ell_1,j_1;\ell_2,j_2}^{j_3,n,m}=-(-1)^{\ell_1\ell_2}\sum_{p=m}^{m+n}{(-1)^p\over (p-m)!}\lmd_{\ell_2,j_2;\ell_1,j_1}^{j_3,m+n-p,p}\eqno(4.7)$$
for $\ell_s\in\Bbb{Z}_2,\;j_s\in I_{i_s}$ and $m,n\in\Bbb{N}$. Next we want to find the condition for (1.4). Observe that
\begin{eqnarray*}& &Y^+(\vs_{\ell_1,j_1},z_1)Y^+(\vs_{\ell_2,j_2},z_2)\vs_{\ell_3,j_3}\\&=&
Y^+(\vs_{\ell_1,j_1},z_1)(\sum_{j_4\in I_{\ell_2+\ell_3},\;m_2,n_2\in\Bbb{N}}\lmd_{\ell_2,j_2;\ell_3,j_3}^{j_4,n_2,m_2}\ptl^{n_2}\vs_{\ell_2+\ell_3,j_4}z_2^{-m_2-1})\\&=&\sum_{j_4\in I_{\ell_2+\ell_3},\;m_2,n_2\in\Bbb{N}}\lmd_{\ell_2,j_2;\ell_3,j_3}^{j_4,n_2,m_2}Y^+(\vs_{\ell_1,j_1},z_1)\ptl^{n_2}\vs_{\ell_2+\ell_3,j_4}z_2^{-m_2-1}\hspace{3cm}\end{eqnarray*}
\begin{eqnarray*}&=&\sum_{j_4\in I_{\ell_2+\ell_3},\;p,m_2,n_2\in\Bbb{N}}\lmd_{\ell_2,j_2;\ell_3,j_3}^{j_4,n_2,m_2}(^{n_2}_{\;p})\ptl^{n_2-p}(-\ptl_{z_1})^pY^+(\vs_{\ell_1,j_1},z_1)\vs_{\ell_2+\ell_3,j_4}z_2^{-m_2-1}\\&=&\sum_{j_4\in I_{\ell_2+\ell_3},\;j_5\in I_{\ell_1+\ell_2+\ell_3},\;p,m_1,m_2,n_1,n_2\in\Bbb{N}}(^{n_2}_{\;p})(m_1+1)(m_1+2)\cdots (m_1+p)\\& & \lmd_{\ell_2,j_2;\ell_3,j_3}^{j_4,n_2,m_2}\lmd_{\ell_1,j_1;\ell_2+\ell_3,j_4}^{j_5,n_1,m_1}\ptl^{n_1+n_2-p}\vs_{\ell_1+\ell_2+\ell_3,j_5}z_1^{-m_1-p-1}z_2^{-m_2-1}\\&=&\sum_{j_4\in I_{\ell_2+\ell_3},\;j_5\in I_{\ell_1+\ell_2+\ell_3},\;p,m_1,n_1,m_2,n_2\in\Bbb{N}}(^{n_2-n_1+p}_{\;\;\;\;\;\;p})m_1(m_1-1)\cdots (m_1-p+1)\\& & \lmd_{\ell_2,j_2;\ell_3,j_3}^{j_4,n_2-n_1+p,m_2}\lmd_{\ell_1,j_1;\ell_2+\ell_3,j_4}^{j_5,n_1,m_1-p}\ptl^{n_2}\vs_{\ell_1+\ell_2+\ell_3,j_5}z_1^{-m_1-1}z_2^{-m_2-1}.\hspace{4.7cm}(4.8)\end{eqnarray*}
Exchanging indices 1 and 2 in (4.8) and then $n_1$ and $n_2$, we have
\begin{eqnarray*}& &Y^+(\vs_{\ell_2,j_2},z_2)Y^+(\vs_{\ell_1,j_1},z_1)\vs_{\ell_3,j_3}\\&=&\sum_{j_4\in I_{\ell_1+\ell_3},\;j_5\in I_{\ell_1+\ell_2+\ell_3},\;p,m_1,m_2,n_1,n_2\in\Bbb{N}}(^{n_1-n_2+p}_{\;\;\;\;\;\;\;p})m_2(m_2-1)\cdots (m_2-p+1)\\& & \lmd_{\ell_1,j_1;\ell_3,j_3}^{j_4,n_1-n_2+p,m_1}\lmd_{\ell_2,j_2;\ell_1+\ell_3,j_4}^{j_5,n_2,m_2-p}\ptl^{n_1}\vs_{\ell_1+\ell_2+\ell_3,j_5}z_1^{-m_1-1}z_2^{-m_2-1}\\&=&\sum_{j_4\in I_{\ell_1+\ell_3},\;j_5\in I_{\ell_1+\ell_2+\ell_3},\;p,m_1,m_2,n_1,n_2\in\Bbb{N}}(^{n_2-n_1+p}_{\;\;\;\;\;\;\;p})m_2(m_2-1)\cdots (m_2-p+1)\\& &\lmd_{\ell_1,j_1;\ell_3,j_3}^{j_4,n_2-n_1+p,m_1}\lmd_{\ell_2,j_2;\ell_1+\ell_3,j_4}^{j_5,n_1,m_2-p}\ptl^{n_2}\vs_{\ell_1+\ell_2+\ell_3,j_5}z_1^{-m_1-1}z_2^{-m_2-1},\hspace{4.6cm}(4.9)\end{eqnarray*}
Moreover, we have
\begin{eqnarray*}& &\rd_x{Y^+(Y^+(\vs_{\ell_1,j_1},z_1-x)\vs_{\ell_2,j_2},x)\vs_{\ell_3,j_3}\over z_2-x}\\&=&\rd_x{1\over z_2-x}(\sum_{j_4\in I_{\ell_1+\ell_2},\;m_1,n_1\in\Bbb{N}}\lmd_{\ell_1,j_1;\ell_2,j_2}^{j_4,n_1,m_1}Y^+(\ptl^{n_1}\vs_{\ell_1+\ell_2,j_4},x)\vs_{\ell_3,j_3}(z_1-x)^{-m_1-1})\\&=&
\rd_x{1\over z_2-x}(\sum_{j_4\in I_{\ell_1+\ell_2},\;m_1,n_1\in\Bbb{N}}\lmd_{\ell_1,j_1;\ell_2,j_2}^{j_4,n_1,m_1}(d^{n_1}Y^+(\vs_{\ell_1+\ell_2,j_4},x)\vs_{\ell_3,j_3}/dx^{n_1})(z_1-x)^{-m_1-1})\\&=&
\rd_x{1\over z_2-x}(\sum_{j_4\in I_{\ell_1+\ell_2},\;j_5\in I_{\ell_1+\ell_2+\ell_3},\;m_1,m_2,n_2,n_1\in\Bbb{N}}\lmd_{\ell_1,j_1;\ell_2,j_2}^{j_4,n_1,m_1}\lmd_{\ell_1+\ell_2,j_4;\ell_3,j_3}^{j_5,n_2,m_2}\\ & &(-1)^{n_1}(m_2+1)(m_2+2)\cdots (m_2+n_1)\ptl^{n_2}\vs_{\ell_1+\ell_2+\ell_3,j_5}x^{-m_2
-n_1-1}(z_1-x)^{-m_1-1})\\&=&
\sum_{j_4\in I_{\ell_1+\ell_2},\;j_5\in I_{\ell_1+\ell_2+\ell_3},\;m_1,m_2,n_2,n_1\in\Bbb{N}}\lmd_{\ell_1,j_1;\ell_2,j_2}^{j_4,n_1,m_1}\lmd_{\ell_1+\ell_2,j_4;\ell_3,j_3}^{j_5,n_2,m_2}\\& &(-1)^{n_1}(m_2+1)(m_2+2)\cdots (m_2+n_1)\ptl^{n_2}\vs_{\ell_1+\ell_2+\ell_3,j_5}\sum_{p=0}^{m_2+n_1}(^{m_1+p}_{\;\;\;\;\:p})z_1^{-m_1-p-1}z_2^{p-m_2-n_1-1}
\\&=&
\sum_{j_4\in I_{\ell_1+\ell_2},\;j_5\in I_{\ell_1+\ell_2+\ell_3},\;p,m_1,m_2,n_2,n_1\in\Bbb{N}}(-1)^{n_1}(m_2+p)\cdots (m_2+p-n_1+1)(^{m_1}_{\;p})\\ & &
\lmd_{\ell_1,j_1;\ell_2,j_2}^{j_4,n_1,m_1-p}\lmd_{\ell_1+\ell_2,j_4;\ell_3,j_3}^{j_5,n_2,m_2+p-n_1}
\ptl^{n_2}\vs_{\ell_1+\ell_2+\ell_3,j_5}z_1^{-m_1-1}z_2^{-m_2-1}.\hspace{4.8cm}(4.10)\end{eqnarray*}
Thus (1.4) is equivalent to:
\begin{eqnarray*}& &
\sum_{j_4\in I_{\ell_2+\ell_3},\;p,n_1\in\Bbb{N}}(^{n_2-n_1+p}_{\;\;\;\;\;\;\;p})m_1(m_1-1)\cdots (m_1-p+1) \lmd_{\ell_2,j_2;\ell_3,j_3}^{j_4,n_2-n_1+p,m_2}\lmd_{\ell_1,j_1;\ell_2+\ell_3,j_4}^{j_5,n_1,m_1-p}\\& &-(-1)^{\ell_1\ell_2}
\sum_{j_4\in I_{\ell_1+\ell_3},\;p,n_1\in\Bbb{N}}(^{n_2-n_1+p}_{\;\;\;\;\;\;\;p})m_2(m_2-1)\cdots (m_2-p+1)\\& &\lmd_{\ell_1,j_1;\ell_3,j_3}^{j_4,n_2-n_1+p,m_1}\lmd_{\ell_2,j_2;\ell_1+\ell_3,j_4}^{j_5,n_1,m_2-p}\\ &=&
\sum_{j_4\in I_{\ell_1+\ell_2},\;p,n_1\in\Bbb{N}}(-1)^{n_1}(m_2+p)\cdots (m_2+p-n_1+1)(^{m_1}_{\;p})\\& &\lmd_{\ell_1,j_1;\ell_2,j_2}^{j_4,n_1,m_1-p}\lmd_{\ell_1+\ell_2,j_4;\ell_3,j_3}^{j_5,n_2,m_2+p-n_1}\hspace{9.7cm}(4.11)\end{eqnarray*}
for $\ell_s\in\Bbb{Z}_2,\;j_s\in I_{\ell_s}$ and $m_1,m_2,n_2\in\Bbb{N}$.

Next we shall use the settings in Section 2. Recall that the definition of the exterior algebra ${\cal E}$ in (2.1) and (2.2), and $\{\psi_{i,j}\mid j\in I_i\}$ is a set of $C^{\infty}$-functions in a real variable $x$ with the ranges in ${\cal E}_i$, for $i\in\Bbb{Z}_2$. We defined $\psi^{(n)}_{i,j}$ in (2.5) and defined  ${\cal A}$ to be  the associative subalgebra of the algebra of functions in the real variable $x$ with the range in ${\cal E}$ generated by 
$\{\psi^{(n)}_{i,j}\mid n\in\Bbb{N},\;i\in\Bbb{Z}_2,\;j\in I_i\}$ with the $\Bbb{Z}_2$-grading in (2.7).  Moreover, we defined $\bar{\cal G}$ in (2.16)-(2.18) and its Lie bracket was given by (2.23).  We set $\td{\cal A}={\cal A}/(d/dx)({\cal A})$ and defined an action of $\bar{\cal G}$ on $\td{\cal A}$ in (2.42). The space $\Omega$ was defined in (2.45) and (2.46).

We associate the above data $(R,V,Y^+(\cdot,z))$ with a matrix differential operator $H:\Omega_0\rightarrow \bar{\cal G}_0$  as follows: for $\vec{v}\in\Omega_0$,
$$H(\vec{v})_{\ell_1,j_1}=\sum_{\ell_2\in\Bbb{Z}_2,\;j_2\in I_{\ell_2},\;j_3\in I_{\ell_1+\ell_2},\;m,n\in\Bbb{N}}{1\over n!}\lmd_{\ell_2,j_2;\ell_1,j_1}^{j_3,m,n}\psi_{\ell_1+\ell_2,j_3}^{(m)}\left({d\over dx}\right)^n(v_{\ell_2,j_2})\eqno(4.12)$$
for $\ell_1\in\Bbb{Z}_2$ and $j_1\in I_{\ell_1}$.
By (2.55), the skew-symmetry (2.54) of $H$ is equivalent to
\begin{eqnarray*}& &\sum_{j_3\in I_{\ell_1+\ell_2},\;m,n\in\Bbb{N}}{1\over m!}\lmd_{\ell_2,j_2;\ell_1,j_1}^{j_3,n,m}\psi_{\ell_1+\ell_2,j_3}^{(n)}\left({d\over dx}\right)^m\\&=&-(-1)^{\ell_1\ell_2}
\sum_{j_3\in I_{\ell_1+\ell_2},\;m,n\in\Bbb{N}}{1\over m!}\lmd_{\ell_1,j_1;\ell_2,j_2}^{j_3,n,m}\left(-{d\over dx}\right)^m\psi_{\ell_1+\ell_2,j_3}^{(n)}\\&=&-(-1)^{\ell_1\ell_2}
\sum_{j_3\in I_{\ell_1+\ell_2},\;m,n\in\Bbb{N}}{1\over m!}\lmd_{\ell_1,j_1;\ell_2,j_2}^{j_3,n,m}(-1)^m\sum_{p=0}^m(^m_{\;p})\psi_{\ell_1+\ell_2,j_3}^{(m+n-p)}\left({d\over dx}\right)^p\\&=&-(-1)^{\ell_1\ell_2}
\sum_{j_3\in I_{\ell_1+\ell_2},\;m,n,p\in\Bbb{N}}{(-1)^p\over p!}(^{\;p}_m)\lmd_{\ell_1,j_1;\ell_2,j_2}^{j_3,n,p}\psi_{\ell_1+\ell_2,j_3}^{(n+p-m)}\left({d\over dx}\right)^m\\&=&-(-1)^{\ell_1\ell_2}\sum_{j_3\in I_{\ell_1+\ell_2},\;m,n,p\in\Bbb{N}}{(-1)^p\over p!}(^{\;p}_m)\lmd_{\ell_1,j_1;\ell_2,j_2}^{j_3,m+n-p,p}\psi_{\ell_1+\ell_2,j_3}^{(n)}\left({d\over dx}\right)^m\hspace{5cm}\end{eqnarray*}
\begin{eqnarray*}
&=&-(-1)^{\ell_1\ell_2}\sum_{j_3\in I_{\ell_1+\ell_2},\;m,n,p\in\Bbb{N},\;p\geq m}{(-1)^p\over m!(p-m)!}\lmd_{\ell_1,j_1;\ell_2,j_2}^{j_3,m+n-p,p}\psi_{\ell_1+\ell_2,j_3}^{(n)}\left({d\over dx}\right)^m,\hspace{1.4cm}(4.13)\end{eqnarray*}
which is equivalent to (4.7). Thus the skew-symmetry of $H$ is equivalent to (1.3).

For convenience, we denote
$$u^{(n)}=\left\{\begin{array}{ll}(d/dx)^n(u)&\mbox{if}\;\;n\in\Bbb{N},\\0&\mbox{otherwise}\end{array}\right.\qquad\for\;\;u\in{\cal A}.\eqno(4.14)$$
Let $\vec{u}_p=\{u_{i,j}^p\}\in \Omega_0$ with $p\in\ol{1,3}$. Using (4.6) and (4.14),  we have
\begin{eqnarray*}& &\vec{u}_3(\ptl_{H(\vec{u}_1)}(H)(\vec{u}_2))\\&=&\sum_{\ell_2,\ell_3\in\Bbb{Z}_2,\;j_2\in I_{\ell_2},\;j_3\in I_{\ell_3},\;j_4\in I_{\ell_2+\ell_3},\;m_2,n_2\in\Bbb{N}}{1\over m_2!}\lmd_{\ell_2,j_2;\ell_3,j_3}^{j_4,n_2,m_2}\{\ptl_{H(\vec{u}_1)}(\psi_{\ell_2+\ell_3,j_4}^{(n_2)})(u_{\ell_2,j_2}^2)^{(m_2)}u^3_{\ell_3,j_3}\}^{\sim}
\\&=&\sum_{\ell_p\in \Bbb{Z}_2,\;j_p\in I_{\ell_p},\;p\in\ol{1,3};\;j_4\in I_{\ell_2+\ell_3},\;j_5\in I_{\ell_1+\ell_2+\ell_3};\;m_1,m_2,n_1,n_2\in \Bbb{N}}{1\over m_1!m_2!}\\& &\lmd_{\ell_2,j_2;\ell_3,j_3}^{j_4,n_2,m_2}\lmd_{\ell_1,j_1;\ell_2+\ell_3,j_4}^{j_5,n_1,m_1}\{[\psi_{\ell_1+\ell_2+\ell_3,j_5}^{(n_1)}(u_{\ell_1,j_1}^1)^{(m_1)}]^{(n_2)}(u_{\ell_2,j_2}^2)^{(m_2)}u^3_{\ell_3,j_3}\}^{\sim}\\&=&\sum_{\ell_p\in \Bbb{Z}_2,\;j_p\in I_{\ell_p},\;m_1,m_2,n_p\in \Bbb{N};\;p\in\ol{1,3};\;j_4\in I_{\ell_2+\ell_3},\;j_5\in I_{\ell_1+\ell_2+\ell_3}}{1\over m_1!m_2!}(^{n_2}_{n_3})\\& &\lmd_{\ell_2,j_2;\ell_3,j_3}^{j_4,n_2,m_2}\lmd_{\ell_1,j_1;\ell_2+\ell_3,j_4}^{j_5,n_1,m_1}\{\psi_{\ell_1+\ell_2+\ell_3,j_5}^{(n_1+n_3)}(u_{\ell_1,j_1}^1)^{(m_1+n_2-n_3)}(u_{\ell_2,j_2}^2)^{(m_2)}u^3_{\ell_3,j_3}\}^{\sim}\hspace{4cm}\end{eqnarray*}
\begin{eqnarray*}&=&\sum_{\ell_p\in \Bbb{Z}_2,\;j_p\in I_{\ell_p},\;m_1,m_2,n_p\in \Bbb{N};\;p\in\ol{1,3};\;j_4\in I_{\ell_2+\ell_3},\;j_5\in I_{\ell_1+\ell_2+\ell_3}}{1\over m_1!m_2!}(^{\;\;\;n_2}_{n_3-n_1})\\& &\lmd_{\ell_2,j_2;\ell_3,j_3}^{j_4,n_2,m_2}\lmd_{\ell_1,j_1;\ell_2+\ell_3,j_4}^{j_5,n_1,m_1}\{\psi_{\ell_1+\ell_2+\ell_3,j_5}^{(n_3)}(u_{\ell_1,j_1}^1)^{(m_1+n_1+n_2-n_3)}(u_{\ell_2,j_2}^2)^{(m_2)}u^3_{\ell_3,j_3}\}^{\sim}
\\&=&\sum_{\ell_p\in \Bbb{Z}_2,\;j_p\in I_{\ell_p},\;m_1,m_2,n_p\in \Bbb{N};\;p\in\ol{1,3};\;j_4\in I_{\ell_2+\ell_3},\;j_5\in I_{\ell_1+\ell_2+\ell_3}}{1\over m_2!(m_1+n_3-n_1-n_2)!}(^{\;\;\;n_2}_{n_3-n_1})\\& &\lmd_{\ell_2,j_2;\ell_3,j_3}^{j_4,n_2,m_2}\lmd_{\ell_1,j_1;\ell_2+\ell_3,j_4}^{j_5,n_1,m_1+n_3-n_1-n_2}\{\psi_{\ell_1+\ell_2+\ell_3,j_5}^{(n_3)}(u_{\ell_1,j_1}^1)^{(m_1)}(u_{\ell_2,j_2}^2)^{(m_2)}u^3_{\ell_3,j_3}\}^{\sim}\hspace{2.1cm}(4.15)\end{eqnarray*}
by (2.19), (2.56) and (4.12). Similarly, we have
\begin{eqnarray*}& &\vec{u}_1(\ptl_{H(\vec{u}_2)}(H)(\vec{u}_3))\\&=&-\vec{u}_3(\ptl_{H(\vec{u}_2)}(H)(\vec{u}_1))
\\&=&\sum_{\ell_p\in \Bbb{Z}_2,\;j_p\in I_{\ell_p},\;m_1,m_2,n_p\in \Bbb{N};\;p\in\ol{1,3};\;j_4\in I_{\ell_1+\ell_3},\;j_5\in I_{\ell_1+\ell_2+\ell_3}}{-1\over (m_2+n_3-n_2-n_1)!m_1!}(^{\;\;\;n_1}_{n_3-n_2})\\& &\lmd_{\ell_1,j_1;\ell_3,j_3}^{j_4,n_1,m_1}\lmd_{\ell_2,j_2;\ell_1+\ell_3,j_4}^{j_5,n_2,m_2+n_3-n_1-n_2}\{\psi_{\ell_1+\ell_2+\ell_3,j_5}^{(n_3)}(u_{\ell_2,j_2}^2)^{(m_2)}(u_{\ell_1,j_1}^1)^{(m_1)}u^3_{\ell_3,j_3}\}^{\sim}\\&=&\sum_{\ell_p\in \Bbb{Z}_2,\;j_p\in I_{\ell_p},\;m_1,m_2,n_p\in \Bbb{N};\;p\in\ol{1,3};\;j_4\in I_{\ell_1+\ell_3},\;j_5\in I_{\ell_1+\ell_2+\ell_3}}{(-1)^{\ell_1\ell_2+1}\over m_1!(m_2+n_3-n_2-n_1)!}(^{\;\;\;n_1}_{n_3-n_2})\\& &\lmd_{\ell_1,j_1;\ell_3,j_3}^{j_4,n_1,m_1}\lmd_{\ell_2,j_2;\ell_1+\ell_3,j_4}^{j_5,n_2,m_2+n_3-n_1-n_2}\{\psi_{\ell_1+\ell_2+\ell_3,j_5}^{(n_3)}(u_{\ell_1,j_1}^1)^{(m_1)}(u_{\ell_2,j_2}^2)^{(m_2)}u^3_{\ell_3,j_3}\}^{\sim}\hspace{2.1cm}(4.16)\end{eqnarray*}
by (2.18), (2.48) and (2.67). Now we assume that $H$ is skew-symmetric. 
 Furthermore,
\begin{eqnarray*}& &\vec{u}_2(\ptl_{H(\vec{u}_3)}(H)(\vec{u}_1))\\&=&\sum_{\ell_1,\ell_2\in\Bbb{Z}_2,\;j_1\in I_{\ell_1},\;j_2\in I_{\ell_2},\;j_4\in I_{\ell_1+\ell_2},\;m_1,n_1\in\Bbb{N}}{1\over m_1!}\\& &\lmd_{\ell_1,j_1;\ell_2,j_2}^{j_4,n_1,m_1}\{
\ptl_{H(\vec{u}_3)}(\psi_{\ell_1+\ell_2,j_4}^{(n_1)})(u_{\ell_1,j_1}^1)^{(m_1)}u^2_{\ell_2,j_2}\}^{\sim}\\&=&\sum_{\ell_1,\ell_2\in\Bbb{Z}_2,\;j_1\in I_{\ell_1},\;j_2\in I_{\ell_2},\;j_4\in I_{\ell_1+\ell_2},\;m_1,n_1\in\Bbb{N}}{1\over m_1!}\\& &\lmd_{\ell_1,j_1;\ell_2,j_2}^{j_4,n_1,m_1}\{
(H(\vec{u}_3)_{\ell_1+\ell_3,j_4})^{(n_1)}(u_{\ell_1,j_1}^1)^{(m_1)}u^2_{\ell_2,j_2}\}^{\sim}
\\&=&\sum_{\ell_1,\ell_2\in\Bbb{Z}_2,\;j_1\in I_{\ell_1},\;j_2\in I_{\ell_2},\;j_4\in I_{\ell_1+\ell_2},\;m_1,n_1\in\Bbb{N}}{(-1)^{n_1}\over m_1!}\\& &\lmd_{\ell_1,j_1;\ell_2,j_2}^{j_4,n_1,m_1}\{
H(\vec{u}_3)_{\ell_1+\ell_3,j_4}[(u_{\ell_1,j_1}^1)^{(m_1)}u^2_{\ell_2,j_2}]^{(n_1)}\}^{\sim}\\&=&\sum_{\ell_p\in\Bbb{Z}_2,\;j_p\in I_{\ell_p};\;p\in\ol{1,3};\;j_4\in I_{\ell_1+\ell_2},\;m_1,m_2,n_1,n_2\in\Bbb{N}}{(-1)^{n_1+1}\over m_1!m_2!}\\& &\lmd_{\ell_1,j_1;\ell_2,j_2}^{j_4,n_1,m_1}
\lmd_{\ell_1+\ell_2,j_4;\ell_3,j_3}^{j_5,n_2,m_2}\{\psi_{\ell_1+\ell_2+\ell_3,j_5}^{(n_2)}[(u_{\ell_1,j_1}^1)^{(m_1)}u^2_{\ell_2,j_2}]^{(m_2+n_1)}u^3_{\ell_3,j_3}\}^{\sim}\\&=&\sum_{\ell_p\in\Bbb{Z}_2,\;j_p\in I_{\ell_p},\;m_1,m_2,n_p\in \Bbb{N};\;p\in\ol{1,3};\;j_4\in I_{\ell_1+\ell_2}}{(-1)^{n_1+1}\over m_1!m_2!}(^{m_2+n_1}_{\;\;\;\;n_3})\\& &\lmd_{\ell_1,j_1;\ell_2,j_2}^{j_4,n_1,m_1}
\lmd_{\ell_1+\ell_2,j_4;\ell_3,j_3}^{j_5,n_2,m_2}\{\psi_{\ell_1+\ell_2+\ell_3,j_5}^{(n_2)}(u_{\ell_1,j_1}^1)^{(m_1+m_2+n_1-n_3)}(u^2_{\ell_2,j_2})^{(n_3)}u^3_{\ell_3,j_3}\}^{\sim}\hspace{4cm}\end{eqnarray*}
\begin{eqnarray*}&=&\sum_{\ell_p\in\Bbb{Z}_2,\;j_p\in I_{\ell_p},\;m_1,m_2,n_p\in \Bbb{N};\;p\in\ol{1,3};\;j_4\in I_{\ell_1+\ell_2}}{(-1)^{n_1+1}\over m_1!n_3!}(^{n_1+n_3}_{\;\;\;m_2})\\& &\lmd_{\ell_1,j_1;\ell_2,j_2}^{j_4,n_1,m_1}
\lmd_{\ell_1+\ell_2,j_4;\ell_3,j_3}^{j_5,n_2,n_3}\psi_{\ell_1+\ell_2+\ell_3,j_5}^{(n_2)}(u_{\ell_1,j_1}^1)^{(m_1+n_1+n_3-m_2)}(u^2_{\ell_2,j_2})^{(m_2)}u^3_{\ell_3,j_3}\}^{\sim}
\\&=&\sum_{\ell_p\in\Bbb{Z}_2,\;j_p\in I_{\ell_p},\;m_1,m_2,n_p\in \Bbb{N};\;p\in\ol{1,3};\;j_4\in I_{\ell_1+\ell_2}}{(-1)^{n_1+1}\over m_1!n_2!}(^{n_1+n_2}_{\;\;\;m_2})\\& &\lmd_{\ell_1,j_1;\ell_2,j_2}^{j_4,n_1,m_1}
\lmd_{\ell_1+\ell_2,j_4;\ell_3,j_3}^{j_5,n_3,n_2}\{\psi_{\ell_1+\ell_2+\ell_3,j_5}^{(n_3)}(u_{\ell_1,j_1}^1)^{(m_1+n_1+n_2-m_2)}(u^2_{\ell_2,j_2})^{(m_2)}u^3_{\ell_3,j_3}\}^{\sim}
\\&=&\sum_{\ell_p\in\Bbb{Z}_2,\;j_p\in I_{\ell_p},\;m_1,m_2,n_p\in \Bbb{N};\;p\in\ol{1,3};\;j_4\in I_{\ell_1+\ell_2}}{(-1)^{n_1+1}\over n_2!(m_1+m_2-n_1-n_2)!}(^{n_1+n_2}_{\;\;\;m_2})\\& &\lmd_{\ell_1,j_1;\ell_2,j_2}^{j_4,n_1,m_1+m_2-n_1-n_2}
\lmd_{\ell_1+\ell_2,j_4;\ell_3,j_3}^{j_5,n_3,n_2}\{\psi_{\ell_1+\ell_2+\ell_3,j_5}^{(n_3)}(u_{\ell_1,j_1}^1)^{(m_1)}(u^2_{\ell_2,j_2})^{(m_2)}u^3_{\ell_3,j_3}\}^{\sim}.\hspace{1.4cm}(4.17)\end{eqnarray*}
by (2.18), (2.43), (2.48) and (2.54). So (2.62) is equivalent to
\begin{eqnarray*}& &\sum_{j_4\in I_{\ell_2+\ell_3},\;n_1,n_2\in \Bbb{N}}{1\over m_2!(m_1+n_3-n_1-n_2)!}(^{\;\;\;n_2}_{n_3-n_1})\lmd_{\ell_2,j_2;\ell_3,j_3}^{j_4,n_2,m_2}\lmd_{\ell_1,j_1;\ell_2+\ell_3,j_4}^{j_5,n_1,m_1+n_3-n_1-n_2}\\& &
-(-1)^{\ell_1\ell_2}\sum_{j_4\in I_{\ell_1+\ell_3},\;n_1,n_2\in \Bbb{N}}{1\over m_1!(m_2+n_3-n_2-n_1)!}(^{\;\;\;n_1}_{n_3-n_2})\lmd_{\ell_1,j_1;\ell_3,j_3}^{j_4,n_1,m_1}\lmd_{\ell_2,j_2;\ell_1+\ell_3,j_4}^{j_5,n_2,m_2+n_3-n_1-n_2}
\hspace{5cm}\end{eqnarray*}
\begin{eqnarray*}&=&\sum_{j_4\in I_{\ell_1+\ell_2},\;n_1,n_2\in\Bbb{N}}{(-1)^{n_1}\over n_2!(m_1+m_2-n_1-n_2)!}(^{n_1+n_2}_{\;\;\;m_2})\\& &\lmd_{\ell_1,j_1;\ell_2,j_2}^{j_4,n_1,m_1+m_2-n_1-n_2}
\lmd_{\ell_1+\ell_2,j_4;\ell_3,j_3}^{j_5,n_3,n_2}\hspace{8.6cm}(4.18)\end{eqnarray*}
for $\ell_p\in\Bbb{Z}_2,\;j_p\in I_p$ and $m_1,m_2,n_3\in\Bbb{N}$.

On the other hand, we can do some changes of indices and combinatorics on (4.11) as follows.
Changing $n_2$ to $n_3$ in (4.11), we have
\begin{eqnarray*}& &
\sum_{j_4\in I_{\ell_2+\ell_3},\;p,n_1\in\Bbb{N}}(^{n_3-n_1+p}_{\;\;\;\;\;\;\;p})m_1(m_1-1)\cdots (m_1-p+1) \lmd_{\ell_2,j_2;\ell_3,j_3}^{j_4,n_3-n_1+p,m_2}\lmd_{\ell_1,j_1;\ell_2+\ell_3,j_4}^{j_5,n_1,m_1-p}\\& &-(-1)^{\ell_1\ell_2}
\sum_{j_4\in I_{\ell_1+\ell_3},\;p,n_1\in\Bbb{N}}(^{n_3-n_1+p}_{\;\;\;\;\;\;\;p})m_2(m_2-1)\cdots (m_2-p+1)\\& &\lmd_{\ell_1,j_1;\ell_3,j_3}^{j_4,n_3-n_1+p,m_1}\lmd_{\ell_2,j_2;\ell_1+\ell_3,j_4}^{j_5,n_1,m_2-p}\\&=&
\sum_{j_4\in I_{\ell_1+\ell_2},\;p,n_1\in\Bbb{N}}(-1)^{n_1}(m_2+p)\cdots (m_2+p-n_1+1)\\& &(^{m_1}_{\;p})\lmd_{\ell_1,j_1;\ell_2,j_2}^{j_4,n_1,m_1-p}\lmd_{\ell_1+\ell_2,j_4;\ell_3,j_3}^{j_5,n_3,m_2+p-n_1}.\hspace{8.9cm}(4.19)\end{eqnarray*}
Then we change $p$ to $n_1+n_2-n_3$ in the first two summations and $p$ to $n_1+n_2-m_2$ in third summation and obtain 
\begin{eqnarray*}& &
\sum_{j_4\in I_{\ell_2+\ell_3},\;n_1,n_2\in\Bbb{N}}(^{\;\;\;\;\;\;\;n_2}_{n_1+n_2-n_3})m_1(m_1-1)\cdots (m_1+n_3-n_1-n_2+1)\hspace{7cm}\end{eqnarray*}
\begin{eqnarray*}& & \lmd_{\ell_2,j_2;\ell_3,j_3}^{j_4,n_2,m_2}\lmd_{\ell_1,j_1;\ell_2+\ell_3,j_4}^{j_5,n_1,m_1+n_3-n_1-n_2}-(-1)^{\ell_1\ell_2}\sum_{j_4\in I_{\ell_1+\ell_3},\;n_1,n_2\in\Bbb{N}}(^{\;\;\;\;\;\;\;n_2}_{n_1+n_2-n_3})\\& &m_2(m_2-1)\cdots (m_2+n_3-n_1-n_2+1)\lmd_{\ell_1,j_1;\ell_3,j_3}^{j_4,n_2,m_1}\lmd_{\ell_2,j_2;\ell_1+\ell_3,j_4}^{j_5,n_1,m_2+n_3-n_1-n_2}\\&=&
\sum_{j_4\in I_{\ell_1+\ell_2},\;n_1,n_2\in\Bbb{N}}(-1)^{n_1}(n_1+n_2)\cdots (n_2+1)(^{\;\;\;\;\;\;\;m_1}_{n_1+n_2-m_2})\lmd_{\ell_1,j_1;\ell_2,j_2}^{j_4,n_1,m_1+m_2-n_1-n_2}\\&&\lmd_{\ell_1+\ell_2,j_4;\ell_3,j_3}^{j_5,n_3,n_2}.\hspace{11.7cm}(4.20)\end{eqnarray*}
Moreover, we exchange $n_1$ and $n_2$ in the second summation and get
\begin{eqnarray*}& &
\sum_{j_4\in I_{\ell_2+\ell_3},\;n_1,n_2\in\Bbb{N}}(^{\;\;\;\;\;\;\;n_2}_{n_1+n_2-n_3})m_1(m_1-1)\cdots (m_1+n_3-n_1-n_2+1) \\& &\lmd_{\ell_2,j_2;\ell_3,j_3}^{j_4,n_2,m_2}\lmd_{\ell_1,j_1;\ell_2+\ell_3,j_4}^{j_5,n_1,m_1+n_3-n_1-n_2}-(-1)^{\ell_1\ell_2}
\sum_{j_4\in I_{\ell_1+\ell_3},\;n_1,n_2\in\Bbb{N}}(^{\;\;\;\;\;\;\;n_1}_{n_1+n_2-n_3})\\& &m_2(m_2-1)\cdots (m_2+n_3-n_1-n_2+1)\lmd_{\ell_1,j_1;\ell_3,j_3}^{j_4,n_1,m_1}\lmd_{\ell_2,j_2;\ell_1+\ell_3,j_4}^{j_5,n_2,m_2+n_3-n_1-n_2}\\&=&
\sum_{j_4\in I_{\ell_1+\ell_2},\;n_1,n_2\in\Bbb{N}}(-1)^{n_1}(n_1+n_2)\cdots (n_2+1)(^{\;\;\;\;\;\;\;m_1}_{n_1+n_2-m_2})\lmd_{\ell_1,j_1;\ell_2,j_2}^{j_4,n_1,m_1+m_2-n_1-n_2}\\& &\lmd_{\ell_1+\ell_2,j_4;\ell_3,j_3}^{j_5,n_3,n_2}.\hspace{11.7cm}(4.21)\end{eqnarray*}
Furthermore, we rewrite (4.21) as
\begin{eqnarray*}& &
\sum_{j_4\in I_{\ell_2+\ell_3},\;n_1,n_2\in\Bbb{N}}(^{\;\;\;n_2}_{n_3-n_1}){m_1!\over(m_1+n_3-n_1-n_2)!} \lmd_{\ell_2,j_2;\ell_3,j_3}^{j_4,n_2,m_2}\lmd_{\ell_1,j_1;\ell_2+\ell_3,j_4}^{j_5,n_1,m_1+n_3-n_1-n_2}\\& &-(-1)^{\ell_1\ell_2}
\sum_{j_4\in I_{\ell_1+\ell_3},\;n_1,n_2\in\Bbb{N}}(^{\;\;\;n_1}_{n_3-n_2}){m_2!\over(m_2+n_3-n_1-n_2)!}\lmd_{\ell_1,j_1;\ell_3,j_3}^{j_4,n_1,m_1}\lmd_{\ell_2,j_2;\ell_1+\ell_3,j_4}^{j_5,n_2,m_2+n_3-n_1-n_2}\\&=&
\sum_{j_4\in I_{\ell_1+\ell_2},\;n_1,n_2\in\Bbb{N}}(-1)^{n_1}{(n_1+n_2)!m_1!\over n_2!(m_1+m_2-n_1-n_2)!(n_1+n_2-m_2)!}\\& &\lmd_{\ell_1,j_1;\ell_2,j_2}^{j_4,n_1,m_1+m_2-n_1-n_2}\lmd_{\ell_1+\ell_2,j_4;\ell_3,j_3}^{j_5,n_3,n_2}.\hspace{8.3cm}(4.22)\end{eqnarray*}
Dividing (4.22) by $m_1!m_2!$, we obtain 
\begin{eqnarray*}& &
\sum_{j_4\in I_{\ell_2+\ell_3},\;n_1,n_2\in\Bbb{N}}{1\over m_2!(m_1+n_3-n_1-n_2)!}(^{\;\;\;n_2}_{n_3-n_1}) \lmd_{\ell_2,j_2;\ell_3,j_3}^{j_4,n_2,m_2}\lmd_{\ell_1,j_1;\ell_2+\ell_3,j_4}^{j_5,n_1,m_1+n_3-n_1-n_2}\\& &-
\sum_{j_4\in I_{\ell_1+\ell_3},\;n_1,n_2\in\Bbb{N}}{(-1)^{\ell_1\ell_2}\over m_1!(m_2+n_3-n_1-n_2)!}(^{\;\;\;n_1}_{n_3-n_2})\lmd_{\ell_1,j_1;\ell_3,j_3}^{j_4,n_1,m_1}\lmd_{\ell_2,j_2;\ell_1+\ell_3,j_4}^{j_5,n_2,m_2+n_3-n_1-n_2}\\&=&
\sum_{j_4\in I_{\ell_1+\ell_2},\;n_1,n_2\in\Bbb{N}}{(-1)^{n_1}(n_1+n_2)!\over n_2!(m_1+m_2-n_1-n_2)!m_2!(n_1+n_2-m_2)!}\\& &\lmd_{\ell_1,j_1;\ell_2,j_2}^{j_4,n_1,m_1+m_2-n_1-n_2}\lmd_{\ell_1+\ell_2,j_4;\ell_3,j_3}^{j_5,n_3,n_2}\\&=&\sum_{j_4\in I_{\ell_1+\ell_2},\;n_1,n_2\in\Bbb{N}}{(-1)^{n_1}\over n_2!(m_1+m_2-n_1-n_2)!}(^{n_1+n_2}_{\;\;\;m_2})\\& &\lmd_{\ell_1,j_1;\ell_2,j_2}^{j_4,n_1,m_1+m_2-n_1-n_2}
\lmd_{\ell_1+\ell_2,j_4;\ell_3,j_3}^{j_5,n_3,n_2},\hspace{8.4cm}(4.23)\end{eqnarray*}
which is (4.18). Thus (1.4) is equivalent to (2.62) when (1.3) holds or equivalently $H$ in (4.12) is skew-symmetric.
We summarize the above results as our main theorem in this section:
\psp

{\bf Theorem 4.1}. {\it Let} $R$ {\it be a} $\Bbb{Z}_2$-{\it graded free} $\Bbb{C}[\ptl]$-{\it module over its} $\Bbb{Z}_2$-{\it graded subspace} $V$ {\it and let} $Y^+(\cdot,z):V\rightarrow LM(V,R[z^{-1}])$ {\it be any given linear map satisfying (4.1). We extend} $Y^+(\cdot,z)$ {\it to a linear map} $Y^+(\cdot,z):R\rightarrow LM(V,R[z^{-1}])$ {\it by (4.2) and (4.3), and define a matrix differential operator} $H$ {\it by (4.12). Then the family} $(R,\ptl,Y^+(\cdot,z))$ {\it forms a conformal superalgebras if and only if} $H$ {\it is a Hamiltonian superoperator}.
\psp

{\bf Remark 4.2}. (a) In [GDo], Gel'fand and Dorfman classified a certain type of Hamiltonian operator by introducing a certain algebraic structure. Balinskii and Novikov [BN] determined certain Poisson brackets of hydrodynamic type by the same algebraic structure. This coincidence is essentially a special example of our correspondence between the Hamiltonian superoperator $H$ of the form (4.12) and the conformal superalgebra $(R,\ptl,Y^+(\cdot,z))$ determined by (4.2)-(4.4). If $R$ is not a free $\Bbb{C}[\ptl]$-module over any subspace, we can establish the analogous correspondence by introducing Hamiltonian superoperators associated with certain quotient modules of ${\cal A}$, whose kernels naturally correspond systems ordinary differential equations. 
 
(b) There is a special case in which we can still establish a direct correspondence.
An element $v\in R_0$ of a conformal superalgebra $(R,\ptl,Y^+(\cdot,z))$ is called a {\it central element} if
$$\ptl v=0,\;\;Y^+(u,z)v=0\qquad\for\;\;u\in R.\eqno(4.24)$$
For a central element $v$, we also have
$$Y^+(v,z)u=0\qquad\for\;\;u\in R\eqno(4.25)$$
by (1.3). Let $(R,\ptl,Y^+(\cdot,z))$ be a conforma superalgebra such that
$$R=\Bbb{C}[\ptl]V\oplus \Bbb{C}{\bf 1},\eqno(4.26)$$
where $V$ is a $\Bbb{Z}_2$-graded subspace of $R$, $\Bbb{C}[\ptl]V$ is a free $\Bbb{C}[\ptl]$-module and ${\bf 1}$ is a central element.
Let $\{\vs_{i,j}\in I_i\}$ be a fixed basis of $V_i$ for $i\in \Bbb{Z}_2$, where $I_1$ and $I_2$ are index sets. We write
$$Y^+(\vs_{\ell_1,j_1},z)\vs_{\ell_2,j_2}=\sum_{j_3\in I_{\ell_1+\ell_2},\;m\in\Bbb{N}}(\mu_{\ell_1,j_1;\ell_2,j_2}^m{\bf 1}+\sum_{n=0}^{\infty}\lmd_{\ell_1,j_1;\ell_2,j_2}^{j_3,n,m}\ptl^n\vs_{\ell_1+\ell_2,j_3})z^{-m-1}
\eqno(4.27)$$
for $\ell_1,\ell_2\in \Bbb{Z}_2$ and $j_p\in I_{\ell_p}$, where $\mu_{\ell_1,j_1;\ell_2,j_2}^m,\lmd_{\ell_1,j_1;\ell_2,j_2}^{j_3,n,m}\in\Bbb{C}$.
We define a matrix differential operator $H:\Omega_0\rightarrow \bar{\cal G}_0$  as follows: for $\vec{v}\in\Omega_0$,
\begin{eqnarray*}& &H(\vec{v})_{\ell_1,j_1}=\sum_{\ell_2\in\Bbb{Z}_2,\;j_2\in I_{\ell_2},\;j_3\in I_{\ell_1+\ell_2},n\in\Bbb{N}}{1\over n!}[\mu_{\ell_2,j_2;\ell_1,j_1}^n\\& &+\sum_{m=0}^{\infty}\lmd_{\ell_2,j_2;\ell_1,j_1}^{j_3,m,n}\psi_{\ell_1+\ell_2,j_3}^{(m)}]\left({d\over dx}\right)^n(v_{\ell_2,j_2})\hspace{7.4cm}(4.28)\end{eqnarray*}
for $\ell_1\in\Bbb{Z}_2$ and $j_1\in I_{\ell_1}$. Then it can be proved that $H$ is a Hamiltonian superoperator. Conversely, a Hamiltonian superoperator $H$ of the form (4.28) determined a conformal superalgebras through (4.24), (4.26) and (4.27). 

(c) Let ${\cal G}$ be a Lie algebras with an  invariant symmetric bilinear form $\la\cdot,\cdot\ra$. Let $R_1=\{0\},\;V={\cal G}$ in (4.26) and the map $Y^+(\cdot,z)$ is determined by
$$Y^+(u,z)v=[u,v]z^{-1}+\la u,v\ra{\bf 1} z^{-2}\qquad\for\;\;u,v\in {\cal G}.\eqno(4.29)$$
This is a conformal algebra generating affine Lie algebra. Let $\{\vs_j\in I\}$ be a fixed basis of ${\cal G}$ and write
$$[\vs_{j_1},\vs_{j_2}]=\sum_{j_3\in I}\lmd_{j_1,j_2}^{j_3}\vs_{j_3},\;\;\mu_{j_1,j_2}=\la \vs_{j_1},\vs_{j_2}\ra\qquad\for\;\;j_1,j_2\in I.\eqno(4.30)$$
Let $I_1=\emptyset$ and $I=I_0$ in Section 2. We denote $\psi_{0,j}=\psi_j$ for $j\in I$. Then the corresponding Hamiltonian operator $H$ is given by
$$H(\vec{u})_{j_1}=\sum_{j_2\in I}(\lmd_{j_2,j_1}^{j_3}\psi_{j_3}u_{j_2}+\mu_{j_2,j_1}(d/dx)(u_{j_2}))\eqno(4.31)$$
for $\vec{u}\in\Omega_0=\Omega$.

\vspace{1cm}

\noindent{\Large \bf References}

\hspace{0.5cm}

\begin{description}

\item[{[BVV]}] B. Bakalov, V. G. Kac and A. Voronov, Cohomology of conformal algebras, {\it Commun. Math. Phys.} {\bf 200} (1999), 561-598.

\item[{[BN]}] A. A. Balinskii and S. P. Novikov, Poisson brackets of hydrodynamic type, Frobenius algebras and Lie algebras, {\it Soviet Math. Dokl.} Vol. {\bf 32} (1985), No. {\bf 1}, 228-231.

\item[{[Bo]}] R. E. Borcherds, Vertex algebras, Kac-Moody algebras, and the Monster,
{\it Proc. Natl. Acad. Sci. USA} {\bf 83} (1986), 3068-3071.

\item[{[CK]}] S.-J. Cheng and V. G. Kac, A new $N=6$ superconformal algebras, {\it Commun. Math. Phys.} {\bf 186} (1997), 219-231.

\item[{[Da1]}]
Yu. L. Daletsky, Lie superalgebras in Hamiltonian operator theory, In: {\it Nonlinear and Turbulent Processes in Physics}, ed. V. E. Zakharov, 19984, pp. 1307-1312.

\item[{[Da2]}]---, Hamiltonian operators in graded formal calculus of variables,
{\it Func. Anal. Appl.} {\bf 20} (1986), 136-138.

\item[{[DK]}] A. D'Andrea and V. G. Kac, Structure theory of finite conformal algebras, {\it Selecta Math (N.S.)} {\bf 4} (1998), 377-418.

\item[{[De]}] B. DeWitt, {\it Supermanifolds,} Second Edition, Cambridge University Press, 1992. 

\item[{[FLM]}] ---, {\it Vertex Operator
Algebras and the Monster}, Pure and Applied Math. Academic Press, 1988.

\item[{[K1]}] V. G. Kac, {\it Vertex algebras for beginners}, University lectures series, Vol {\bf 10}, AMS. Providence RI, 1996.

\item[{[K2]}] ---, Superconformal algebras and transitive group actions on quadrics, {\it Commun. Math. Phys.} {\bf 186} (1997), 233-252.

\item[{[K3]}] ---, Idea of locality, {\it Physical Applications and Mathematical Aspects of Geometry, Groups and Algebras}, Doebener et al eds., World Scientific Publishers, 1997, 16-32.

\item[{[GDi1]}] I. M. Gel'fand and L. A. Dikii, Asymptotic behaviour of the resolvent of Sturm-Liouville equations and the algebra of the Korteweg-de Vries equations, {\it Russian Math. Surveys} {\bf 30:5} (1975), 77-113.

\item[{[GDi2]}] ---, A Lie algebra structure in a formal variational Calculation, {\it Func. Anal. Appl.}  {\bf 10} (1976), 16-22.

\item[{[GDo]}] 
I. M. Gel'fand and I. Ya. Dorfman, Hamiltonian operators and algebraic structures related to them, {\it Funkts. Anal. Prilozhen}  {\bf 13} (1979), 13-30.

\item[{[KT]}] V. G. Kac and I. T. Todorov, Superconformal current algebras and their unitary representations, {\it Commun. Math. Phys}. {\bf 102} (1985), 337-347.

\item[{[M]}] P. Mathieu, Supersymmetry extension of the Korteweg-de Vries equation, {\it J. Math. Phys.} {\bf 29} (11) (1988), 2499-2507.

\item[{[X1]}] X. Xu, Hamiltonian superoperators, {\it J. Phys A: Math. \& Gen.} {\bf 28} No. 6 (1995).

\item[{[X2]}] ---, Variational calculus of supervariables and related algebraic structures, {\it J. Algebra}, in press; preprint was circulated in January 1995.

\item[{[X3]}] ---, {\it Introduction to Vertex Operator Superalgebras and Their Modules}, Kluwer Academic Publishers, Dordrecht/Boston/London, 1998.

 \item[{[X4]}] ---,  Quadratic conformal superalgebras, {\it J. Algebra}, to appear.

\item[{[X5]}] ---, Simple conformal superalgebras of finite growth, {\it submitted}.

\end{description}
\end{document}